\documentclass{elsarticle}

\makeatletter
\def\ps@pprintTitle{
 \let\@oddhead\@empty
 \let\@evenhead\@empty
 \def\@oddfoot{\centerline{\thepage}}%
 \let\@evenfoot\@oddfoot}
\makeatother

\usepackage[utf8]{inputenc}        

\usepackage{amsmath,amssymb,mathrsfs}
\usepackage{bm}                    

\usepackage{algorithm}
\usepackage{algpseudocode}         

\usepackage{booktabs}
\usepackage{multirow}
\usepackage{makecell}
\usepackage{enumerate}

\usepackage{float}
\usepackage{upref}
\usepackage{geometry}            
\usepackage[colorlinks]{hyperref}

\numberwithin{equation}{section}

\DeclareMathOperator{\tr}{tr}
\DeclareMathOperator{\erf}{erf}
\DeclareMathOperator{\erfc}{erfc}

\newcommand{\os}{\hspace{1ex}}  

\def\ccm{Center for Computational Mathematics, Flatiron Institute, Simons Foundation, New York, NY, 10010, USA}
\def\sjtu{School of Mathematical Sciences, Shanghai Jiao Tong University, Shanghai, 200240, China}

\begin{document}

\begin{frontmatter}

\title{Fast Ewald Summation with Prolates for Charged Systems in the NPT Ensemble}

\author[ccm,sjtu]{Jiuyang Liang}
\ead{jliang@flatironinstitute.org}

\author[ccm]{Libin Lu}
\ead{llu@flatironinstitute.org}

\author[ccm]{Shidong Jiang\corref{mycor}}
\ead{sjiang@flatironinstitute.org}

\address[ccm]{\ccm}
\address[sjtu]{\sjtu}

\cortext[mycor]{Corresponding author}

\begin{abstract}
We present an NPT extension of Ewald summation with prolates (ESP), a spectrally accurate and scalable particle-mesh method for molecular dynamics simulations of periodic, charged systems. Building on the recently introduced ESP framework, this work focuses on rigorous and thermodynamically consistent pressure/stress evaluation in the isothermal--isobaric ensemble. ESP employs prolate spheroidal wave functions as both splitting and spreading kernels, reducing the Fourier grid size needed to reach a prescribed pressure accuracy compared with current widely used mesh-Ewald methods based on Gaussian splitting and B-spline spreading. We derive a unified pressure-tensor formulation applicable to isotropic, semi-isotropic, anisotropic, and fully flexible cells, and show that the long-range pressure can be evaluated with a single forward FFT followed by diagonal scaling, whereas force evaluation requires both forward and inverse transforms. We provide production implementations in LAMMPS and GROMACS and validate pressure and force accuracy on bulk water, LiTFSI ionic liquids, and a transmembrane system. Benchmarks on up to $3\times 10^3$ CPU cores demonstrate strong scaling and reduced communication cost at matched accuracy, particularly for NPT pressure evaluation.
\end{abstract}


\end{frontmatter}

\section{Introduction}
Molecular dynamics (MD) simulation has become a versatile tool for studying physical, chemical, biological, and materials systems at the atomic scale~\cite{karplus1990molecular,hollingsworth2018molecular}. A central goal in MD is the accurate sampling of statistical ensembles, in which macroscopic states are specified by fixed thermodynamic variables~\cite{frenkel2001understanding}. Many laboratory conditions are well described by constant particle number ($N$), pressure ($P$), and temperature ($T$), corresponding to the isothermal--isobaric (NPT) ensemble. Consequently, NPT simulations are widely used for molecular crystals and biomolecular systems such as ribosomes and transmembrane proteins~\cite{cordova2023atomic,xia1997crystal,bock2013energy}.

Sampling the NPT ensemble relies on thermostats and barostats to control temperature and pressure. A broad range of schemes has been proposed, including deterministic methods such as Nos\'e--Hoover~\cite{nose1983constant} and MTK~\cite{martyna1994constant}, and stochastic methods such as the Langevin piston~\cite{feller1995constant}, COMPEL~\cite{di2015stochastic}, cell-rescaling~\cite{bernetti2020pressure}, and the second-order Langevin sampler preserving positive volume~\cite{li2024second}. Unlike NVE or NVT dynamics, NPT equations of motion depend explicitly on the instantaneous pressure, which must be evaluated at each step of time integration. For systems with only short-range interactions, the pressure tensor can be computed efficiently using the standard virial expression~\cite{frenkel2001understanding}, together with nearest-image corrections under periodic boundary conditions~\cite{Thompson2009JCP}. For Coulomb interactions, however, pressure evaluation is substantially more challenging: the long-range nature of electrostatics makes direct virial calculations scale as $O(N^2)$~\cite{Thompson2009JCP,louwerse2006calculation}.

Classical Ewald summation addresses this difficulty by splitting electrostatics into real-space and Fourier-space contributions~\cite{ewald1921berechnung}. In the NPT setting, pressure can be obtained by differentiating the Helmholtz free energy with respect to volume~\cite{brown1995general,hummer1998pressure}. With lattice-based techniques and fast Fourier transforms (FFTs), mesh Ewald methods achieve $O(N\log N)$ complexity for energies and forces, and can also support pressure evaluation~\cite{di2015stochastic,sega2016pressure}. However, on modern parallel architectures these approaches require communication-heavy all-to-all operations over grid data, which can limit scalability; in large-scale NPT simulations, long-range electrostatics therefore often dominates the overall cost. Compared with the extensive literature on fast methods for energies and forces~\cite{greengard1987fast,greengard1988rapid,Hockney1988Computer,barnes1986nature}, fewer algorithms have been developed that specifically target fast and scalable evaluation of the long-range pressure (or stress) tensor. Recent work to reduce communication cost includes the random batch Ewald method~\cite{liang2022jcp,Liang_2022}, which replaces FFTs with stochastic sampling over a small mini-batch. While efficient in an ensemble-averaged sense, its accuracy is statistical and can be insufficient for pressure-sensitive settings, particularly when cell shape fluctuations are important. An accurate pressure formulation for Coulomb interactions, together with a fast and scalable solver, thus remains a practical challenge for NPT simulations.

In this paper, we extend Ewald summation with prolates (ESP) to MD simulations in the NPT ensemble, with a focus on consistent and efficient evaluation of the instantaneous pressure tensor. ESP differs from traditional mesh Ewald approaches in its use of prolate spheroidal wave functions (PSWFs): PSWFs serve as the kernel for splitting (in place of Gaussian screening) and as the particle-to-grid spreading kernel. The use of PSWFs for kernel splitting was first noted in the dual-space multilevel kernel-splitting (DMK) framework~\cite{jiang2025dual}. We subsequently developed ESP, using PSWFs for both splitting and spreading, for spectrally accurate energy and force evaluation in NVT simulations~\cite{liang2025acceleratingfastewaldsummation}. In the present work, we extend ESP to the NPT ensemble by deriving an NPT-consistent decomposition of the pressure tensor into real-space and Fourier-space contributions. The real-space term decays rapidly and is evaluated by direct truncation, while the Fourier-space term is computed with a 3D spectral solver in which the long-range pressure requires only a single forward FFT followed by diagonal scaling, whereas potentials and forces require both forward and inverse transforms~\cite{liang2025acceleratingfastewaldsummation}. With PSWF-based spreading and a near-optimal support size, ESP reduces the number of grid points per particle and the Fourier grid size needed to reach a prescribed pressure accuracy, providing a unified and efficient framework for energy, force, and pressure evaluation.

We assess the resulting NPT-capable ESP method through systematic numerical experiments and parallel benchmarks (Table~\ref{tbl:expts}). ESP retains the $O(N\log N)$ complexity of mesh-based Ewald methods while reducing both computation and communication at matched pressure accuracy, and we implement it in LAMMPS and GROMACS with only minor modifications to the existing PPPM/PME workflow. As summarized in Table~\ref{tbl:expts}, strong-scaling tests on up to $3\times 10^{3}$ CPU cores show that, for a $2.4\times 10^{7}$-particle water system, ESP is up to three times faster than GROMACS-PME at $4\times 10^{-4}$ accuracy, while at higher accuracy around $10^{-5}$ it achieves $5$--$8\times$ speedups relative to LAMMPS-PPPM. At matched accuracy, ESP typically reduces the required number of Fourier grid points by roughly an order of magnitude. We further validate ESP on bulk SPC/E water, a LiTFSI ionic liquid, and the transmembrane bovine bc$_1$ complex, observing accurate thermodynamic and structural statistics over time scales from femtoseconds to microseconds. Together, these results indicate that ESP is a practical and efficient option for large-scale NPT simulations of charged systems.

\begin{table}
  \centering
  \scriptsize
  \resizebox{0.95\textwidth}{!}{%
  \begin{tabular}{|llcllccll|}
    \hline
    system & $N$ atoms & $\Delta$ & code & cores &
    \multicolumn{1}{c}{PME/PPPM} &
    \multicolumn{1}{c}{ESP} &
    speed-up & Fig.\\
    &&&&& $h$(nm), $P$ & $h$(nm), $P$ &&\\
    \hline
    \multicolumn{9}{|l|}{\rule{0pt}{3ex}\hspace{-2ex}{\bf Strong-scaling tests for large MD:}}\\
    water & 24M & $4\times 10^{-4}$ & G & 96--3k & 0.12, \os 5 & 0.26,\os 5 & 2--3$\times$ & 6b \\
    water & 11M & $2\times 10^{-5}$ & L & 96--3k & 0.067, \os 5 & 0.2,\os 6 & 5--8$\times$ & 7b \\
    \hline
    \multicolumn{9}{|l|}{\rule{0pt}{3ex}\hspace{-2ex}{\bf Long-time MD, standard accuracy:}}\\
    Transmembrane       & 1M & $4\times 10^{-4}$ & G & 960 & 0.12, \os 5 & 0.26, \os 5 & 2$\times$   & 9 \\
    LiTFSI ionic liquid & 1M & $2\times 10^{-4}$ & L & 960 & 0.11, \os 5 & 0.24, \os 5 & 2.9$\times$ & 8 \\
    \hline
  \end{tabular}}
  \caption{Benchmark comparison of the ESP method against native PME and PPPM electrostatics in \textsc{GROMACS} and \textsc{LAMMPS} for the NPT ensemble. Results are grouped by study type
  (strong-scaling for large systems; long-time MD at standard accuracy). $\Delta$ denotes the requested error
  tolerance; for each case, the PME/PPPM grid spacing $h$ and interpolation
  order $P$ are tuned to achieve a comparable error level. Code: G = \textsc{GROMACS}, L = \textsc{LAMMPS}.
  Reported speed-ups are for the pure Coulomb calculation in the strong-scaling tests
  (Figs.~6b, 7b) and for end-to-end MD throughput in the long-time simulations (Figs.~8, 9).}
  \label{tbl:expts}
\end{table}

The remainder of this paper is organized as follows. Section~2 introduces the PSWF-based kernel-splitting formulation for Coulomb interactions. Section~3 derives pressure-tensor expressions consistent with isotropic, semi-isotropic, anisotropic, and fully flexible pressure-coupling schemes. Section~4 describes the resulting ESP algorithm and implementation details. Section~5 presents numerical results assessing accuracy and performance. Section~6 concludes the paper.

\section{Electrostatics and the microscopic pressure}\label{sec::PSWFSplitting}
This section summarizes preliminaries on electrostatic interactions, kernel decompositions, and the evaluation of instantaneous pressure under several common parameterizations of the simulation cell. The corresponding partition functions and their dependence on cell variables are discussed in \ref{app::DistributionFunction}.

\subsection{Coulomb interactions and kernel decomposition}\label{subsec::electrostatic}
Consider a charge-neutral system of $N$ particles at positions $\bm{r}_i\in\mathbb{R}^3$ with charges $q_i$, $i=1,\ldots,N$, in an orthorhombic cell $\Omega$ with side lengths $L_x$, $L_y$, and $L_z$. We impose three-dimensional periodic boundary conditions. Charge neutrality means $\sum_{i=1}^{N} q_i = 0$. The electrostatic potential at particle $i$ is given by the lattice sum
\begin{equation}\label{eq::1.1}
\Phi(\bm{r}_i)=\sum_{j=1}^{N}\sum_{\bm{n}\in\mathbb{Z}^3}{}^{\prime}\frac{q_j}{\left|\bm{r}_{ij}+\bm{n}\circ\bm{L}\right|},
\end{equation}
where $\bm{r}_{ij}:=\bm{r}_i-\bm{r}_j$, $\bm{L}=(L_x,L_y,L_z)$, the prime indicates that the term $i=j$ and $\bm{n}=\bm{0}$ is omitted, and ``$\circ$'' denotes the Hadamard (componentwise) product. The total electrostatic energy is
\begin{equation}
U=\frac{1}{2}\sum_{i=1}^{N}q_i\,\Phi(\bm{r}_i),
\end{equation}
where the factor $1/2$ avoids double counting. The force on particle $i$ is the negative gradient of the energy, $\bm{F}(\bm{r}_i)=-\nabla_{\bm{r}_i}U$.

The lattice sum in \eqref{eq::1.1} is conditionally convergent (its value depends on the summation order), so naive truncation is not reliable. Moreover, the Coulomb kernel is singular at the origin, which complicates a direct Fourier-space treatment.

Classical Ewald summation~\cite{ewald1921berechnung} addresses these issues by splitting the Coulomb kernel into short- and long-range parts,
\begin{equation}
\frac{1}{r}=\frac{\erfc(\sqrt{\alpha}\,r)}{r}+\frac{\erf(\sqrt{\alpha}\,r)}{r}
:=\mathcal{N}(r)+\mathcal{F}(r),
\end{equation}
where
\begin{equation}
\erf(r)=\frac{2}{\sqrt{\pi}}\int_{0}^{r}e^{-t^2}\,dt
\end{equation}
is the error function, $\erfc(r)=1-\erf(r)$ is its complement, and $\alpha>0$ is the splitting parameter. The near-field term $\mathcal{N}(r)$ decays rapidly and can be truncated at a cutoff radius $r_c$, while the smooth far-field term $\mathcal{F}(r)$ is handled in Fourier space. When combined with particle--mesh discretization and FFTs, the long-range component can be evaluated in $O(N\log N)$ time per step. This is the basis of standard mesh-Ewald methods such as PPPM and PME, which are widely used in MD packages including LAMMPS~\cite{plimpton1995fast} and GROMACS~\cite{berendsen1995gromacs}.

As an alternative, the $u$-series method~\cite{predescu2020u} replaces the Ewald split by a sum-of-Gaussians (SOG) approximation. The kernel is represented by a Gaussian series, and near- and far-field contributions are grouped into terms with small and large bandwidths, respectively. This approach can match Ewald accuracy while reducing the Fourier-space cost by roughly a factor of two, and it has motivated follow-up algorithms~\cite{Liang2023SISRBSOG,gao_fast_2025}, theoretical analyses~\cite{liang2025error}, and applications including plasma simulations~\cite{chen2025random} and machine-learning potentials~\cite{ji2025machine}.

\subsection{Instantaneous pressure}\label{subsec::inspress}
In NPT simulations, the time integrator requires both forces and the instantaneous pressure (or pressure tensor), not forces alone. We briefly outline pressure evaluation for periodic systems under fixed and variable cell shapes.

We describe the simulation cell by the cell matrix
\begin{equation}
\bm{h}=[\bm{h}_1,\bm{h}_2,\bm{h}_3]\in\mathbb{R}^{3\times 3},
\end{equation}
whose columns $\bm{h}_j$ span the cell. It is convenient to factor
\begin{equation}
\bm{h}=V^{1/3}\bm{h}_0,\qquad \det(\bm{h}_0)=1,
\end{equation}
so that $V=\det(\bm{h})$ is the volume and $\bm{h}_0$ encodes the shape. Let $\{\bm{r}_1,\ldots,\bm{r}_N\}\equiv\bm{r}_{\mathrm{tot}}$ and $\{\bm{p}_1,\ldots,\bm{p}_N\}\equiv\bm{p}_{\mathrm{tot}}$ denote Cartesian positions and momenta. Using reduced coordinates $\bm{s}_i\in\mathbb{R}^3$, we write
\begin{equation}
\bm{r}_i=\bm{h}\bm{s}_i=\sum_{j=1}^{3}s_i^{(j)}\bm{h}_j,
\end{equation}
where $s_i^{(j)}$ is the $j$th component of $\bm{s}_i$.

\noindent \textbf{Isotropic (cubic) coupling.}
Let $k_{\mathrm{B}}$ be Boltzmann's constant, $T$ the temperature, and $E=K+U$ the total energy. For a cubic cell with isotropic coupling, the thermodynamic pressure is~\cite{frenkel2001understanding}
\begin{equation}
P=\frac{1}{\beta}\frac{\partial \log Q(N,V,T)}{\partial V},
\end{equation}
where $\beta=(k_{\mathrm{B}}T)^{-1}$ and
\begin{equation}
Q(N,V,T)=\int_{\Omega\times\mathbb{R}^{3N}}\exp\!\bigl(-\beta E(\bm{r}_{\mathrm{tot}},\bm{p}_{\mathrm{tot}})\bigr)\,
d\bm{r}_{\mathrm{tot}}\,d\bm{p}_{\mathrm{tot}}
\end{equation}
is the canonical partition function. Because $\bm{r}_i$ depends on $V$, it is convenient to change variables to the reduced conjugate pair
\[
\bm{s}_i=V^{-1/3}\bm{r}_i,\qquad \bm{p}_i^{s}=V^{2/3}m_i\dot{\bm{s}}_i,
\]
where $m_i$ is the particle mass. The kinetic energy expressed in these variables is
\begin{equation}
K(V,\bm{p}^s)=V^{-2/3}\sum_{i=1}^{N}\frac{|\bm{p}_i^{s}|^2}{2m_i}.
\end{equation}
Differentiating $\log Q$ with respect to $V$ and using the definition of ensemble averages yields the instantaneous pressure as the negative volume derivative of the energy at fixed reduced coordinates,
\begin{equation}\label{eq::Pins}
P_{\mathrm{ins}}
= -\left.\frac{\partial (K+U)}{\partial V}\right|_{\{\bm{s}_i\}}
= \frac{1}{3V}\left[\sum_{i=1}^{N}\frac{|\bm{p}_i|^2}{m_i}
-3V\left.\frac{\partial U(\{V^{1/3}\bm{s}_i\};V)}{\partial V}\right|_{\{\bm{s}_i\}}\right]_{\bm{s}_i=V^{-1/3}\bm{r}_i}.
\end{equation}
Evaluating the derivative more explicitly gives the familiar virial form plus an additional term,
\begin{equation}\label{eq::ins}
P_{\mathrm{ins}}
=\frac{1}{3V}\sum_{i=1}^{N}\frac{|\bm{p}_i|^2}{m_i}
-\frac{1}{3V}\sum_{i=1}^{N}\bm{r}_i\cdot\nabla_{\bm{r}_i}U(\{\bm{r}_i\};V)
-\left.\frac{\partial U(\{\bm{r}_i\};V)}{\partial V}\right|_{\{\bm{r}_i\}}.
\end{equation}

Under periodic boundary conditions, the last term in \eqref{eq::ins} represents an extra momentum flux associated with the implicit motion of periodic images induced by volume or shape changes~\cite{Thompson2009JCP}. The correction approach in~\cite{Thompson2009JCP} is not convenient for Coulomb interactions when kernel splitting is used. In Section~\ref{sec::ESPPressure}, we address this by a hybrid strategy that applies \eqref{eq::Pins} and \eqref{eq::ins} to the real- and Fourier-space contributions, respectively.

\noindent \textbf{Semi-isotropic coupling.}
For semi-isotropic coupling (commonly used for interfacial or membrane systems), the lateral area $A$ is allowed to fluctuate while remaining isotropic in the plane, and the cell height $L$ fluctuates independently (so $V=AL$). The corresponding instantaneous pressures can be written as
\begin{equation}
P_{\mathrm{ins},A}=-\left.\frac{1}{L}\frac{\partial E}{\partial A}\right|_{\{\bm{s}_i\},L}
=-\frac{1}{L}\left(\frac{\partial K}{\partial A}+\frac{\partial U}{\partial A}\right),
\quad
P_{\mathrm{ins},L}=-\left.\frac{1}{A}\frac{\partial E}{\partial L}\right|_{\{\bm{s}_i\},A}
=-\frac{1}{A}\left(\frac{\partial K}{\partial L}+\frac{\partial U}{\partial L}\right).
\end{equation}

\noindent \textbf{Anisotropic (orthorhombic) coupling.}
For an orthorhombic cell with independently fluctuating side lengths $(L_x,L_y,L_z)$, the instantaneous pressure tensor is diagonal. Using Greek indices $\alpha,\beta,\gamma\in\{1,2,3\}$ (equivalently $\{x,y,z\}$) for Cartesian components, its diagonal entries are
\begin{equation}\label{eq::diagona}
(\bm{P}_{\mathrm{ins}})_{\alpha\alpha}
= -\left.\frac{L_\alpha}{V}\frac{\partial (K+U)}{\partial L_\alpha}\right|_{\{\bm{s}_i\}}
= \frac{1}{V}\left[\sum_{i=1}^{N}\frac{(p_i)_\alpha^{2}}{m_i}
-\left.\frac{\partial U(\{\bm{h}\bm{s}_i\};\bm{h})}{\partial L_\alpha}\right|_{\bm{s}_i\leftarrow \bm{h}^{-1}\bm{r}_i}\,L_\alpha\right].
\end{equation}

\noindent \textbf{Fully flexible cell.}
For a general (possibly non-orthogonal) cell, we use the change of variables $\bm{r}_i=\bm{h}\bm{s}_i$ and $\bm{p}_i=\bm{h}^{-T}\bm{p}_i^s$, whose Jacobian is one. The kinetic energy is
\begin{equation}
K(\bm{h},\bm{p}^s)=\sum_{i=1}^{N}\frac{\left|\bm{h}^{-T}\bm{p}_i^s\right|^2}{2m_i}.
\end{equation}
The instantaneous pressure (stress) tensor can be written as
\begin{equation}\label{eq::instheta}
(\bm{P}_{\mathrm{ins}})_{\alpha\beta}
=\frac{1}{\det(\bm{h})}\left(\frac{\partial(-K-U)}{\partial \bm{h}}\bm{h}^{T}\right)_{\alpha\beta}
=\frac{1}{\det(\bm{h})}\left(\sum_{i=1}^{N}\frac{\bm{p}_i\otimes\bm{p}_i}{m_i}
-\sum_{\gamma}\frac{\partial U(\{\bm{h}\bm{s}_{i}\},\bm{h})}{\partial h_{\alpha\gamma}}\,h_{\beta\gamma}\right).
\end{equation}
Here the derivative $\partial/\partial \bm{h}$ is taken at fixed reduced variables $\{\bm{s}_i,\bm{p}_i^s\}$, where $\bm{r}_i=\bm{h}\bm{s}_i$ and $\bm{p}_i=\bm{h}^{-T}\bm{p}_i^s$. In practice, the reduced variables are obtained from the MD state via $\bm{s}_i=\bm{h}^{-1}\bm{r}_i$ and $\bm{p}_i^s=\bm{h}^{T}\bm{p}_i$. Greek indices $\alpha,\beta,\gamma\in\{1,2,3\}$ (equivalently $\{x,y,z\}$) denote Cartesian components. When $\bm{h}$ is diagonal, \eqref{eq::diagona} corresponds to the diagonal part of \eqref{eq::instheta}. Expanding the derivative of $U$ with respect to $h_{\alpha\gamma}$ gives
\begin{equation}\label{eq::pressuretensor}
(\bm{P}_{\mathrm{ins}})_{\alpha\beta}
=\frac{1}{\det(\bm{h})}\left[\sum_{i=1}^{N}\left(\frac{\bm{p}_i\otimes\bm{p}_i}{m_i}
-\nabla_{\bm{r}_i}U\otimes \bm{r}_i\right)_{\alpha\beta}
-\left.\sum_{\gamma}\frac{\partial U(\{\bm{r}_i\};\bm{h})}{\partial h_{\alpha\gamma}}\right|_{\{\bm{r}_i\}}\,h_{\beta\gamma}\right].
\end{equation}
As in the isotropic case, the final term in \eqref{eq::pressuretensor} arises under periodic boundary conditions and accounts for the implicit dependence of the electrostatic energy on cell shape through the periodic images.

\section{Pressure calculations with prolate spheroidal wave functions}\label{sec::ESPPressure}

In this section, we introduce a PSWF-based kernel splitting for Coulomb interactions and use it to derive instantaneous pressure (and pressure-tensor) formulas under general pressure-coupling schemes. The derivation explicitly accounts for the self-energy, the zero-frequency (mean-field) contribution, and the treatment of non-neutral charge distributions, yielding expressions that are consistent with periodic boundary conditions and amenable to fast Fourier-space evaluation. We also prove that the PSWF splitting is compatible with the virial theorem, ensuring that the resulting pressure formulas are thermodynamically consistent.

\subsection{PSWF splitting for Coulomb interactions}\label{subsec::pswfdecomp}
Let $\psi_0^c(\cdot)$ denote the order-zero prolate spheroidal wave function (PSWF) with parameter $c>0$; its definition and key properties are summarized in \ref{app::pswf}. Fix a real-space cutoff radius $r_c>0$ (hereafter referred to as the \emph{cutoff radius}). Using $\psi_0^c$, we split the Coulomb kernel into a compactly supported near-field term and a smooth far-field term,
\begin{equation}\label{eq::decomposition}
\mathcal{N}^c(r)=\begin{cases}
\dfrac{1-\phi_0^c(r)}{r},\quad& \text{if}~r\leq r_c,\\[1.5em]
0,\quad & \text{if}~r>r_c,
\end{cases}
\qquad
\mathcal{F}^c(r)=\frac{\phi_0^c(r)}{r},
\end{equation}
where
\begin{equation}
\phi_0^c(r):=\frac{1}{C_0}\int_0^{r/r_c}\psi_0^c(x)\,dx,
\qquad
C_0:=\int_{0}^1\psi_0^c(x)\,dx.
\end{equation}
By construction, $\mathcal{N}^c(r)$ retains the Coulomb singularity but is supported only on $[0,r_c]$, whereas $\mathcal{F}^c(r)$ is smooth at the origin. In particular,
\begin{equation}
\lim_{r\rightarrow 0}\mathcal{F}^c(r)=\frac{\psi_{0}^{c}(0)}{C_0r_c},
\end{equation}
so the far-field kernel is well suited for Fourier-space treatment.

By construction, the PSWF splitting in \eqref{eq::decomposition} is exact for all $r>0$, i.e.,
\begin{equation}
\mathcal{N}^c(r)+\mathcal{F}^c(r)\equiv \frac{1}{r}\qquad (r>0).
\end{equation}
Two additional features are particularly useful for fast particle--mesh algorithms. First, the near-field kernel is \emph{compactly supported} with cutoff radius $r_c$: it vanishes identically for $r>r_c$ and, crucially, satisfies
\begin{equation}
\mathcal{N}^c(r_c)=0.
\end{equation}
Thus the real-space pair potential contributed by $\mathcal{N}^c$ is naturally \emph{zero at the cutoff} without any shift or switching. This is in contrast to the classical Ewald near-field term $\erfc(\sqrt{\alpha}\,r)/r$, which is generally nonzero at $r_c$ and therefore requires a consistent potential shift (and corresponding bookkeeping) in practical implementations such as GROMACS when truncated in real space. Second, the far-field kernel $\mathcal{F}^c$ is effectively bandlimited: its Fourier transform is concentrated in a low-frequency region, which reduces the number of Fourier modes required for a prescribed accuracy. For conventional Gaussian screening (as in standard mesh-Ewald methods), achieving an error level $\varepsilon$ typically requires
\begin{equation}\label{eq:KmaxPPPM}
K_{\max} =  2\log(1/\varepsilon).
\end{equation}
Here $K_{\max}$ denotes the maximum Fourier frequency retained in the spectral (FFT-based) evaluation in each dimension. For PSWF splitting, the analogous bandwidth parameter is $c$, and to reach the same small $\varepsilon$ one typically chooses $c$ such that
\begin{equation}\label{eq::c}
K_{\max}= c \approx \log(1/\varepsilon).
\end{equation}
The approximation becomes more accurate in the high-precision regime (small $\varepsilon$). Consequently, at sufficiently high accuracy the required Fourier grid resolution is reduced by roughly a factor of two in each dimension relative to Gaussian-based splittings, leading to substantial savings in three dimensions. These bandwidth advantages carry over directly to pressure-tensor evaluation, where the Fourier-space contribution is computed on the same reciprocal grid. In the next section, we show how these advantages carry over to pressure-tensor evaluation by deriving an NPT-compatible decomposition into real-space and Fourier-space contributions, each treated with a method suited to its decay and regularity.

\subsection{PSWF-based instantaneous pressure calculation}\label{subsec::pswfpressure}
We now apply the PSWF splitting to derive the instantaneous pressure tensor for a general triclinic cell with cell matrix $\bm{h}=[\bm{h}_1,\bm{h}_2,\bm{h}_3]$. The reciprocal basis vectors are
\begin{equation}\label{eq::vector}
\bm{b}_1=2\pi\frac{\bm{h}_2\times\bm{h}_3}{\det(\bm{h})},\qquad
\bm{b}_2=2\pi\frac{\bm{h}_3\times\bm{h}_1}{\det(\bm{h})},\qquad
\bm{b}_3=2\pi\frac{\bm{h}_1\times\bm{h}_2}{\det(\bm{h})},
\end{equation}
so that $\bm{h}_\alpha\cdot \bm{b}_\beta=2\pi \delta_{\alpha\beta}$ for all $\alpha,\beta\in\{1,2,3\}$, where $\delta_{\alpha\beta}$ is the Kronecker delta. Any reciprocal vector can be written as
\begin{equation}
\bm{k}=\sum_{\alpha=1}^{3}m_{\alpha}\bm{b}_{\alpha},\qquad m_{\alpha}\in\mathbb{Z},
\end{equation}
equivalently $\bm{k}=2\pi \bm{h}^{-T}\bm{m}$ with $\bm{m}=[m_1,m_2,m_3]^{T}\in\mathbb{Z}^3$. We use the Fourier pair on the periodic cell $\Omega=\{\bm{h}\bm{s}:\bm{s}\in[0,1)^3\}$,
\begin{equation}
\widehat{f}(\bm{k})=\int_{\Omega}f(\bm{r})e^{-i\bm{k}\cdot\bm{r}}\,d\bm{r},\qquad
f(\bm{r})=\frac{1}{\det(\bm{h})}\sum_{\bm{k}}\widehat{f}(\bm{k})e^{i\bm{k}\cdot\bm{r}}.
\end{equation}

With the PSWF splitting in \eqref{eq::decomposition}, the Coulomb energy separates into short-range, long-range, and self-energy contributions,
\begin{equation}\label{eq::splitting}
U=U_{\mathcal{N}}+U_{\mathcal{F}}+U_{\text{self}}.
\end{equation}
The short-range energy is
\begin{equation}\label{eq::UN}
U_{\mathcal{N}}
=\frac{1}{2}\sum_{\bm{n}\in\mathbb{Z}^3}{}^{\prime}\sum_{i,j=1}^{N}q_iq_j\,
\mathcal{N}^c\!\left(|\bm{r}_{ij}+\bm{h}\bm{n}|\right),
\end{equation}
where $\bm{r}_{ij}:=\bm{r}_{i}-\bm{r}_j$ and the prime indicates that the term $i=j$ with $\bm{n}=\bm{0}$ is excluded. The long-range energy admits the Fourier representation
\begin{equation}\label{eq::UF}
U_{\mathcal{F}}
=\frac{1}{2\det(\bm{h})}\sum_{\bm{k}\neq\bm{0}}\widehat{\mathcal{F}}^c(\bm{k})\left|\rho(\bm{k})\right|^2
+U_{\mathcal{F}}^{0},
\end{equation}
where
\begin{equation}
\widehat{\mathcal{F}}^c(\bm{k})
=\frac{2\pi \lambda_0}{C_0}\,\frac{\psi_0^{c}(|\bm{k}|r_c/c)}{|\bm{k}|^2},
\qquad
\rho(\bm{k}):=\sum_{i=1}^{N}q_ie^{i\bm{k}\cdot\bm{r}_i}
\end{equation}
is the charge structure factor. The term $U_{\mathcal{F}}^{0}$ collects the contribution associated with the $\bm{k}=\bm{0}$ (zero-frequency) mode, whose value depends on the summation order and the choice of macroscopic boundary conditions~\cite{zhonghan2014JCTC,Liang2023SISRBSOG}. Throughout this work we adopt tin-foil (conducting) boundary conditions, under which the $\bm{k}=\bm{0}$ mode is omitted; accordingly, we set $U_{\mathcal{F}}^{0}\equiv 0$. The self-energy term removes the unphysical self-interaction,
\begin{equation}\label{eq::selfen}
U_{\text{self}}
=-\frac{1}{2}\sum_{i=1}^{N}q_i^2\mathcal{F}^c(0)
=-\frac{\psi_0^c(0)}{2C_0r_c}\sum_{i=1}^{N}q_i^2,
\end{equation}
and does not contribute to the pressure since it is independent of $\bm{h}$.

A useful observation for the pressure derivation is that $\bm{r}_{ij}=\bm{h}(\bm{s}_i-\bm{s}_j)$ and
\begin{equation}
\rho(\bm{k})
=\sum_{i=1}^{N}q_i e^{i(2\pi \bm{h}^{-T}\bm{m})\cdot(\bm{h}\bm{s}_i)}
=\sum_{i=1}^{N}q_i e^{i2\pi \bm{m}^{T}\bm{s}_i},
\end{equation}
so $\rho(\bm{k})$ is independent of $\bm{h}$ when expressed in reduced coordinates.

Using the general NPT pressure formulas from Section~\ref{subsec::inspress}, we decompose the instantaneous pressure tensor into kinetic, short-range, and long-range parts,
\begin{equation}
\bm{P}_{\mathrm{ins}}
=\frac{1}{\det(\bm{h})}\sum_{i=1}^{N}\frac{\bm{p}_i\otimes \bm{p}_i}{m_i}
+\bm{P}_{\mathrm{ins},\mathcal{N}}+\bm{P}_{\mathrm{ins},\mathcal{F}},
\end{equation}
where $\bm{P}_{\mathrm{ins},\mathcal{N}}$ and $\bm{P}_{\mathrm{ins},\mathcal{F}}$ arise from $U_{\mathcal{N}}$ and $U_{\mathcal{F}}$, respectively.

\noindent \textbf{Real-space contribution.}
Applying the cell-derivative form \eqref{eq::pressuretensor} to $U_{\mathcal{N}}$ yields the pairwise expression
\begin{equation}\label{eq::PinsN}
\bm{P}_{\mathrm{ins},\mathcal{N}}
=-\frac{1}{2\det(\bm{h})}\sum_{\bm{n}\in\mathbb{Z}^3}{}^{\prime}\sum_{i,j=1}^{N}q_{i}q_{j}\,
\mathcal{F}_{\mathcal{N}}\!\left(|\bm{r}_{ij}+\bm{h}\bm{n}|\right)\,
\frac{(\bm{r}_{ij}+\bm{h}\bm{n})\otimes (\bm{r}_{ij}+\bm{h}\bm{n})}{|\bm{r}_{ij}+\bm{h}\bm{n}|^3},
\end{equation}
where the scalar weight $\mathcal{F}_{\mathcal{N}}(r)$ is defined by
\begin{equation}\label{eq::mathFN}
\mathcal{F}_{\mathcal{N}}(r)
:=1-\frac{1}{C_0}\int_{0}^{r/r_c}\psi_0^c(x)\,dx+\frac{r}{C_0 r_c}\psi_0^{c}(r/r_c).
\end{equation}
This form is convenient in practice because it is consistent with the corresponding real-space force term (\ref{app::force}) and can be evaluated using the same neighbor list.

\noindent \textbf{Fourier-space contribution.}
A direct pairwise form for the long-range pressure is not practical due to $O(N^2)$ scaling, so we differentiate the Fourier representation \eqref{eq::UF}. Using
\begin{equation}
\frac{\partial |\bm{k}|^2}{\partial \bm{h}}\bm{h}^{T}=-2\bm{k}\otimes\bm{k},
\qquad
\frac{\partial}{\partial \bm{h}}\!\left(\frac{1}{\det(\bm{h})}\right)\bm{h}^{T}
=-\frac{1}{\det(\bm{h})}I,
\end{equation}
we obtain
\begin{equation}\label{eq::UFPree}
\bm{P}_{\mathrm{ins},\mathcal{F}}
=\frac{1}{2\det(\bm{h})^2}\sum_{\bm{k}\neq \bm{0}}|\rho(\bm{k})|^2
\left[\widehat{\mathcal{F}}^c(\bm{k})I
+2\frac{\partial \widehat{\mathcal{F}}^c(\bm{k})}{\partial |\bm{k}|^2}\,\bm{k}\otimes\bm{k}\right],
\end{equation}
where $I$ is the $3\times 3$ identity matrix. A direct calculation gives
\begin{equation}
\frac{\partial \widehat{\mathcal{F}}^c(\bm{k})}{\partial |\bm{k}|^2}
=-\frac{\widehat{\mathcal{F}}^c(\bm{k})}{|\bm{k}|^2}
+\frac{\pi \lambda_0}{C_0}\,\frac{|\bm{k}|r_c}{c}\,
\frac{\psi_0^{c}{}^\prime(|\bm{k}|r_c/c)}{|\bm{k}|^4}.
\end{equation}
Substituting this into \eqref{eq::UFPree} yields the explicit long-range pressure tensor
\begin{equation}\label{eq::Far}
\bm{P}_{\mathrm{ins},\mathcal{F}}
=-\frac{1}{2\det(\bm{h})}\sum_{\bm{k}\neq \bm{0}}|\rho(\bm{k})|^2\left[
\widehat{\mathcal{F}}^c(\bm{k})\left(I-\frac{2\bm{k}\otimes\bm{k}}{|\bm{k}|^2}\right)
+\frac{2\pi\lambda_0}{C_0}\frac{|\bm{k}|r_c}{c}\frac{\psi_{0}^{c}{}^{\prime}(|\bm{k}|r_c/c)}{|\bm{k}|^4}\bm{k}\otimes\bm{k}
\right].
\end{equation}
By \ref{app::pswf}, the above Fourier sum can be truncated at $|\bm{k}|\le c/r_c$. Combining all contributions, we obtain an efficient expression for the instantaneous pressure tensor,
\begin{equation}\label{eq::tensoraniso}
\begin{split}
\bm{P}_{\mathrm{ins}}
=&\frac{1}{\det(\bm{h})}\sum_{i=1}^{N}\frac{\bm{p}_i\otimes \bm{p}_i}{m_i}
-\frac{1}{2\det(\bm{h})}\sum_{\bm{n}\in\mathbb{Z}^3}{}^{\prime}\sum_{i,j=1}^{N}q_{i}q_{j}\,
\mathcal{F}_{\mathcal{N}}\!\left(|\bm{r}_{ij}+\bm{h}\bm{n}|\right)\,
\frac{(\bm{r}_{ij}+\bm{h}\bm{n})\otimes (\bm{r}_{ij}+\bm{h}\bm{n})}{|\bm{r}_{ij}+\bm{h}\bm{n}|^3}\\
&+\frac{1}{2\det(\bm{h})^2}\sum_{\bm{k}\neq \bm{0}}|\rho(\bm{k})|^2\left[
\widehat{\mathcal{F}}^c(\bm{k})\left(I-\frac{2\bm{k}\otimes\bm{k}}{|\bm{k}|^2}\right)
+\frac{2\pi\lambda_0}{C_0}\frac{|\bm{k}|r_c}{c}\frac{\psi_{0}^{c}{}^{\prime}(|\bm{k}|r_c/c)}{|\bm{k}|^4}\bm{k}\otimes\bm{k}
\right].
\end{split}
\end{equation}

\noindent \textbf{Orthorhombic specialization.}
For an orthorhombic cell $\bm{h}=\mathrm{diag}(L_x,L_y,L_z)$, the pressure tensor is diagonal. Writing $V=L_xL_yL_z$ and $k_{\alpha}=2\pi m_{\alpha}/L_{\alpha}$ with $\alpha\in\{x,y,z\}$, we have
\begin{equation}\label{eq::rec}
\begin{split}
(\bm{P}_{\mathrm{ins}})_{\alpha\alpha}
=&\frac{1}{V}\sum_{i=1}^{N}\frac{(p_i)_\alpha^{2}}{m_i}
-\frac{1}{2V}\sum_{\bm{n}\in\mathbb{Z}^3}{}^{\prime}\sum_{i,j=1}^{N}q_iq_j\,
\mathcal{F}_{\mathcal{N}}\!\left(|\bm{r}_{ij}+\bm{h}\bm{n}|\right)\,
\frac{(\bm{r}_{ij}+\bm{h}\bm{n})_{\alpha}^2}{|\bm{r}_{ij}+\bm{h}\bm{n}|^3}\\
&+\frac{1}{2V^2}\sum_{\bm{k}\neq \bm{0}}|\rho(\bm{k})|^2\left[
\widehat{\mathcal{F}}^{c}(\bm{k})\left(1-\frac{2k_{\alpha}^2}{|\bm{k}|^2}\right)
+\frac{2\pi\lambda_0}{C_0}\frac{|\bm{k}|r_c}{c}\frac{\psi_{0}^{c}{}^{\prime}(|\bm{k}|r_c/c)}{|\bm{k}|^4}k_{\alpha}^2
\right],
\end{split}
\end{equation}
and $(\bm{P}_{\mathrm{ins}})_{\alpha\beta}=0$ for $\alpha\neq \beta$. For a cubic cell, the instantaneous scalar pressure is one-third of the trace of \eqref{eq::rec}.

\subsection{Verification of the virial theorem}\label{subsec::virial}
The virial theorem~\cite{frenkel2001understanding} relates the (ensemble-averaged) pressure of a system to its kinetic and potential energies. For Coulomb interactions, the potential is homogeneous of degree $-1$ under uniform scaling: for any $\gamma>0$,
\begin{equation}
U(\{\gamma\bm{r}_i\};\gamma L)=\gamma^{-1}U(\{\bm{r}_i\};L).
\end{equation}
As a consequence, the isotropic pressure satisfies
\begin{equation}
P_{\mathrm{iso}}
=\left\langle \frac{1}{3}\sum_{\alpha\in\{x,y,z\}}(\bm{P}_{\mathrm{ins}})_{\alpha \alpha} \right\rangle
= \frac{N}{\beta\,\det(\bm{h})} + \left\langle \frac{U}{3\,\det (\bm{h})} \right\rangle,
\end{equation}
where $\langle\cdot\rangle$ denotes an ensemble average. Equivalently, using \eqref{eq::instheta}, the Coulomb homogeneity implies
\begin{equation}\label{eq::deth}
-\frac{1}{\det(\bm{h})}\tr\!\left[\frac{\partial U}{\partial \bm{h}}\bm{h}^{T}\right]
=\frac{U}{\det(\bm{h})}.
\end{equation}
With the PSWF splitting, \eqref{eq::deth} can be written as
\begin{equation}\label{eq::PSRela}
\tr\!\left[\bm{P}_{\mathcal{N}}+\bm{P}_{\mathcal{F}}\right]
=\frac{1}{\det(\bm{h})}\left(U_{\mathcal{N}}+U_{\mathcal{F}}+U_{\text{self}}\right).
\end{equation}
If a splitting violates \eqref{eq::PSRela}, the resulting pressure statistics can be biased, compromising the reliability of NPT sampling.

We now verify \eqref{eq::PSRela} for the PSWF splitting. A direct calculation yields
\begin{equation}\label{eq::UNS}
\frac{U_{\mathcal{N}}+U_{\text{self}}}{\det(\bm{h})}-\tr\!\left[\bm{P}_{\mathcal{N}}\right]
=-\frac{1}{2\det(\bm{h})}\sum_{i=1}^{N}q_i\,\xi(\bm{r}_i),
\end{equation}
and
\begin{equation}\label{eq::UFS}
\frac{U_{\mathcal{F}}}{\det(\bm{h})}-\tr\!\left[\bm{P}_{\mathcal{F}}\right]
=-\frac{1}{2\det(\bm{h})^{2}}\sum_{\bm{k}\neq \bm{0}}|\rho(\bm{k})|^2
\left[\frac{2\pi\lambda_0}{C_0}\frac{|\bm{k}|r_c}{c}\frac{\psi_{0}^{c}{}^{\prime}(|\bm{k}|r_c/c)}{|\bm{k}|^2}\right],
\end{equation}
where
\begin{equation}
\xi(\bm{r})
:=\sum_{\bm{n}\in\mathbb{Z}^3}\sum_{j=1}^{N}\frac{q_j}{C_0r_c}\,
\psi_{0}^{c}\!\left(\frac{|\bm{r}-\bm{r}_{j}+\bm{h}\bm{n}|}{r_c}\right)
\end{equation}
is radially symmetric. Using the identity in \ref{App:Fourier} and the change of variables $u=r/r_c$, the Fourier transform of $\xi$ can be written as
\begin{equation}\label{eq::derivative}
\begin{split}
\widehat{\xi}(\bm{k})
&=\sum_{j=1}^{N}\frac{q_j}{C_0r_c}\int_{\Omega}\sum_{\bm{n}\in\mathbb{Z}^3}
\psi_{0}^{c}\!\left(\frac{|\bm{r}-\bm{r}_j+\bm{h}\bm{n}|}{r_c}\right)e^{-i\bm{k}\cdot\bm{r}}\,d\bm{r}\\
&=\sum_{j=1}^{N}\frac{q_j}{C_0r_c}\int_{\mathbb{R}^3}
\psi_{0}^{c}\!\left(\frac{|\bm{r}-\bm{r}_j|}{r_c}\right)e^{-i\bm{k}\cdot(\bm{r}-\bm{r}_j)}\,d\bm{r}\;e^{-i\bm{k}\cdot\bm{r}_j}\\
&=\sum_{j=1}^{N}\frac{2\pi q_jr_ce^{-i\bm{k}\cdot\bm{r}_j}}{C_0|\bm{k}|}
\int_{-1}^{1}\psi_{0}^{c}(u)\sin(|\bm{k}| r_c u)\,u\,du.
\end{split}
\end{equation}
Differentiating \eqref{eq::psidefinit} with respect to $x$ and taking $n=0$ gives
\begin{equation}\label{eq::derivative_psi}
\lambda_{0}\psi_{0}^{c}{}^{\prime}(x)
=-c\int_{-1}^{1}\psi_{0}^{c}(u)\sin(cxu)\,u\,du,
\end{equation}
where we used the parity of the integrand. Substituting $x=|\bm{k}| r_c/c$ in \eqref{eq::derivative_psi} and inserting the result into \eqref{eq::derivative} yields
\begin{equation}
\widehat{\xi}(\bm{k})
=-\sum_{j=1}^{N}\frac{2\pi q_j r_c\lambda_0 e^{-i\bm{k}\cdot\bm{r}_j}}{C_0c|\bm{k}|}\,
\psi_{0}^{c}{}^{\prime}\!\left(\frac{|\bm{k}|r_c}{c}\right).
\end{equation}
Therefore, the Fourier expansion of $\xi$ is
\begin{equation}\label{eq::xir}
\xi(\bm{r})
=-\frac{2\pi r_c\lambda_0}{C_0c \det(\bm{h})}\sum_{\bm{k}\neq \bm{0}}
\frac{\psi_{0}^c{}^{\prime}(|\bm{k}|r_c/c)}{|\bm{k}|}\sum_{j=1}^{N}q_je^{-i\bm{k}\cdot(\bm{r}_j-\bm{r})}.
\end{equation}
Substituting \eqref{eq::xir} into \eqref{eq::UNS} gives
\begin{equation}
\begin{split}
\frac{U_{\mathcal{N}}+U_{\text{self}}}{\det(\bm{h})}-\tr\!\left[\bm{P}_{\mathcal{N}}\right]
&=-\frac{\pi r_c \lambda_0}{C_0c\,\det(\bm{h})^{2}}
\sum_{\bm{k}\neq\bm{0}}\frac{\psi_{0}^c{}^{\prime}(|\bm{k}|r_c/c)}{|\bm{k}|}|\rho(\bm{k})|^2\\
&=-\frac{U_{\mathcal{F}}}{\det(\bm{h})}+\tr\!\left[\bm{P}_{\mathcal{F}}\right],
\end{split}
\end{equation}
which, together with \eqref{eq::UFS}, verifies \eqref{eq::PSRela}. This confirms that the PSWF-based formulation is consistent with the virial theorem.

Recent studies have raised questions about pressure evaluation in systems with long-range interactions. One work~\cite{Onegin_2024} suggested that, in Ewald-type methods, the splitting parameter could be treated as volume-dependent and might influence the pressure, whereas other studies~\cite{li2025comment,zhao2025infinite} argued that the pressure is independent of such algorithmic choices. The verification above supports the latter viewpoint: the parameters introduced in ESP need not be treated as volume-dependent in pressure calculations.

\subsection{Correction for non-neutral systems}\label{subsec::correction}
The derivations above assume charge neutrality, $\sum_{i=1}^{N}q_i=0$. If the net charge is nonzero, the $\bm{k}=\bm{0}$ mode renders the Fourier-space energy (and hence the pressure) ill-defined, and an additional convention is required. Let
\begin{equation}
Q_{\mathrm{tot}}:=\sum_{i=1}^{N}q_{i}
\end{equation}
denote the total charge. A standard choice is to introduce a uniform neutralizing background charge density
\[
\rho_{\mathrm{bg}}=-\frac{Q_{\mathrm{tot}}}{V},
\]
so that the combined system is neutral.

With the PSWF splitting, the background-corrected energy can be written as
\begin{equation}\label{eq:Corr}
U_{\mathrm{corr}}
=U_{\mathcal{N}}+U_{\mathcal{F}}+U_{\mathrm{self}}
+U_{\mathrm{c\text{-}b}}^{\mathrm{corr}}+U_{\mathrm{b\text{-}b}}^{\mathrm{corr}}.
\end{equation}
Here $U_{\mathrm{c\text{-}b}}^{\mathrm{corr}}$ and $U_{\mathrm{b\text{-}b}}^{\mathrm{corr}}$ account for charge--background and background--background interactions, respectively, computed using the near-field kernel $\mathcal{N}^c$ (equivalently, the real-space cutoff at $r_c$ is enforced). Specifically,
\begin{equation}
\begin{split}
U_{\mathrm{c\text{-}b}}^{\mathrm{corr}}
&:=\sum_{i=1}^{N}q_i\int_{\Omega}\sum_{\bm{n}\in\mathbb{Z}^3}
\mathcal{N}^c\!\left(|\bm{r}-\bm{r}_i+\bm{h}\bm{n}|\right)\,\rho_{\mathrm{bg}}\,d\bm{r}\\
&=\sum_{i=1}^{N}q_i\rho_{\mathrm{bg}}\int_{\mathbb{R}^3}\mathcal{N}^c(|\bm{r}|)\,d\bm{r}\\
&=-\frac{4\pi Q_{\mathrm{tot}}^2}{V}\left(\frac{r_c^2}{2}-\int_{0}^{r_c}r\,\phi_{0}^{c}(r)\,dr\right),
\end{split}
\end{equation}
and
\begin{equation}
\begin{split}
U_{\mathrm{b\text{-}b}}^{\mathrm{corr}}
&:=\frac{1}{2}\int_{\Omega}\rho_{\mathrm{bg}}\int_{\Omega}\sum_{\bm{n}\in\mathbb{Z}^3}
\mathcal{N}^c\!\left(|\bm{r}-\bm{r}^\prime+\bm{h}\bm{n}|\right)\,\rho_{\mathrm{bg}}\,d\bm{r}^{\prime}d\bm{r}\\
&=\frac{1}{2}\rho_{\mathrm{bg}}^{2}\,V\int_{\mathbb{R}^3}\mathcal{N}^c(|\bm{r}|)\,d\bm{r}\\
&=\frac{2\pi Q_{\mathrm{tot}}^2}{V}\left(\frac{r_c^2}{2}-\int_{0}^{r_c}r\,\phi_0^c(r)\,dr\right).
\end{split}
\end{equation}
The one-dimensional integral can be evaluated accurately with standard quadrature (e.g., Gauss--Legendre). These background terms remove the singular behavior associated with the $\bm{k}=\bm{0}$ mode of the charge--charge interaction, so the Fourier sum can be taken over $\bm{k}\neq\bm{0}$ without an additional zero-mode treatment. Moreover,
\[
\nabla_{\bm{r}_i}\!\left(U_{\mathrm{c\text{-}b}}^{\mathrm{corr}}+U_{\mathrm{b\text{-}b}}^{\mathrm{corr}}\right)\equiv 0,
\]
so the background correction does not modify particle forces.

We use \eqref{eq::instheta} to obtain the corresponding correction to the instantaneous pressure tensor. Differentiating the energy terms gives
\begin{equation}
\bm{P}_{\mathrm{ins}}^{\mathrm{corr}}
=\frac{1}{\det(\bm{h})}\sum_{i=1}^{N}\frac{\bm{p}_i\otimes \bm{p}_i}{m_i}
+\bm{P}_{\mathrm{ins},\mathcal{N}}+\bm{P}_{\mathrm{ins},\mathcal{F}}
+\bm{P}_{\mathrm{c\text{-}b}}^{\mathrm{corr}}+\bm{P}_{\mathrm{b\text{-}b}}^{\mathrm{corr}},
\end{equation}
with
\begin{equation}
\begin{split}
\bm{P}_{\mathrm{c\text{-}b}}^{\mathrm{corr}}
&=-\frac{1}{\det(\bm{h})}\frac{\partial U_{\mathrm{c\text{-}b}}^{\mathrm{corr}}(\bm{h})}{\partial\bm{h}}\bm{h}^{T}\\
&=-\frac{4\pi Q_{\mathrm{tot}}^2}{\det(\bm{h})^{2}}\left(\frac{r_c^2}{2}-\int_{0}^{r_c}r\,\phi_0^{c}(r)\,dr\right)I,
\end{split}
\end{equation}
and
\begin{equation}
\begin{split}
\bm{P}_{\mathrm{b\text{-}b}}^{\mathrm{corr}}
&=-\frac{1}{\det(\bm{h})}\frac{\partial U_{\mathrm{b\text{-}b}}^{\mathrm{corr}}(\bm{h})}{\partial\bm{h}}\bm{h}^{T}\\
&=\frac{2\pi Q_{\mathrm{tot}}^2}{\det(\bm{h})^{2}}\left(\frac{r_c^2}{2}-\int_{0}^{r_c}r\,\phi_0^{c}(r)\,dr\right)I.
\end{split}
\end{equation}
These corrections are spatially homogeneous and contribute only to the isotropic (diagonal) part of the pressure tensor.

\section{Fast algorithm}\label{sec::FastAlgorithm}
In this section, we present an FFT-accelerated method for computing the long-range (Fourier-space) contribution to the instantaneous pressure tensor. While FFT-based pressure formulations are well established~\cite{essmann1995smooth,di2015stochastic,sega2016pressure}, the ESP formulation differs in two main ways. First, we replace the classical Gaussian Ewald decomposition with the PSWF-based splitting introduced in Section~\ref{sec::ESPPressure}. Second, we use PSWFs as the charge-spreading (window) kernel on the uniform mesh, in place of the B-spline windows used in current PME/PPPM implementations. We also summarize implementation details and describe how to extract local (per-particle or per-group) contributions from the global Fourier-space pressure.

\subsection{An Ewald summation with prolates method for pressure calculations}\label{subsec::FFTSpectral}
We introduce a window function $W(\bm{r})$ and its Fourier transform $\widehat{W}(\bm{k})$. Starting from the Fourier-space pressure term $\bm{P}_{\mathrm{ins},\mathcal{F}}$ in \eqref{eq::tensoraniso}, we insert the identity
\begin{equation}
1\equiv \widehat{W}(\bm{k})^{-2}\,\widehat{W}(\bm{k})^{2},
\end{equation}
which yields
\begin{equation}\label{eq::Pressure*}
\bm{P}_{\mathrm{ins},\mathcal{F}}
=\frac{1}{2\det(\bm{h})^{2}}
\sum_{\bm{k}\neq \bm{0}}\widehat{W}(\bm{k})^{-2}\,
\widehat{\bm{\mathcal{P}}}^{c}(\bm{k})
\left|\sum_{j=1}^{N}q_j\widehat{W}(\bm{k})e^{-i\bm{k}\cdot\bm{r}_j}\right|^2,
\end{equation}
where the mode-by-mode pressure kernel $\widehat{\bm{\mathcal{P}}}^{c}(\bm{k})$ is defined by
\begin{equation}
\widehat{\bm{\mathcal{P}}}^{c}(\bm{k})
:=\widehat{\mathcal{F}}^c(\bm{k})\left(I-\frac{2\bm{k}\otimes\bm{k}}{|\bm{k}|^2}\right)
+\frac{2\pi\lambda_0}{C_0}\frac{|\bm{k}|r_c}{c}\frac{\psi_{0}^{c}{}^{\prime}(|\bm{k}|r_c/c)}{|\bm{k}|^4}\bm{k}\otimes\bm{k}.
\end{equation}
This factorization suggests a particle--mesh spectral method: $\widehat{W}(\bm{k})$ is incorporated into the particle-to-grid spreading step, while $\widehat{W}(\bm{k})^{-2}$ appears as a diagonal (mode-by-mode) scaling.

The long-range pressure evaluation consists of four steps.

\textbf{Step 1 (Spreading).} Define the grid charge density
\begin{equation}\label{eq::gridding}
\rho_{\mathrm{grid}}(\bm{r})
=\sum_{j=1}^{N}q_j\,W(\bm{r}-\bm{r}_j)_{*},
\end{equation}
where $(\cdot)_{*}$ denotes periodization of $W$. If the support of $W$ is smaller than half the box length in each direction, only the central image contributes. The quantity $\rho_{\mathrm{grid}}$ is evaluated on a uniform mesh.

\textbf{Step 2 (3D FFT).} Apply a single forward 3D FFT to obtain $\widehat{\rho}_{\mathrm{grid}}(\bm{k})$.

\textbf{Step 3 (Diagonal scaling).} For each Fourier mode, compute
\begin{equation}\label{eq::diag}
\widehat{\bm{\Pi}}(\bm{k})
:=\frac{1}{2\det(\bm{h})^{2}}\,
\widehat{W}(\bm{k})^{-2}\,\widehat{\bm{\mathcal{P}}}^{c}(\bm{k})\,
\left|\widehat{\rho}_{\mathrm{grid}}(\bm{k})\right|^{2}.
\end{equation}

\textbf{Step 4 (Fourier collection).} Sum over modes to obtain the long-range pressure contribution:
\begin{equation}\label{eq::PFfinal}
\bm{P}_{\mathrm{ins},\mathcal{F}}=\sum_{\bm{k}\neq\bm{0}}\widehat{\bm{\Pi}}(\bm{k}).
\end{equation}

A key point is that this four-step procedure requires only one forward FFT. After the diagonal scaling in Step~3, the pressure tensor is obtained directly by a Fourier-space reduction, and no inverse FFT is needed. If $W$ is smooth and compactly supported, the method converges spectrally with respect to grid resolution. The overall procedure is summarized in Algorithm~\ref{al::Fourier}. In ESP, we choose a PSWF-based window $W$; its definition and implementation details are given in Sections~\ref{subsec::PSWFwindow} and \ref{subsec::implement}, respectively.

\begin{algorithm}[H]
\caption{Four-step ESP method for long-range pressure}\label{al::Fourier}
\begin{algorithmic}[1]
\State (Spreading) Evaluate $\rho_{\mathrm{grid}}$ on a uniform mesh using \eqref{eq::gridding}.
\State (3D FFT) Compute $\widehat{\rho}_{\mathrm{grid}}(\bm{k})$ by a forward 3D FFT.
\State (Diagonal scaling) Compute $\widehat{\bm{\Pi}}(\bm{k})$ mode-by-mode using \eqref{eq::diag}.
\State (Fourier collection) Reduce over $\bm{k}\neq\bm{0}$ using \eqref{eq::PFfinal}.
\end{algorithmic}
\end{algorithm}

\subsection{The PSWF window function}\label{subsec::PSWFwindow}
In particle--mesh Ewald methods, the window function controls both accuracy and cost through (i) the real-space support used in spreading/interpolation and (ii) the decay of its Fourier transform, which determines aliasing errors. Common choices include Gaussian windows~\cite{lindbo2011spectral}, B-splines~\cite{Darden1993JCP}, Kaiser--Bessel (KB) windows~\cite{kaiser1980use}, and the ``exponential of semicircle'' (ES) window~\cite{barnett_parallel_2019}. In ESP, we adopt PSWFs as the window function. Their provably optimal spectral concentration among compactly supported functions under a bandlimit constraint enables smaller spreading support at a prescribed accuracy, thereby reducing particle--mesh coupling and communication.

Consider a uniform Cartesian grid on the primary cell $\Omega$, with $M_d$ subintervals in each direction $d\in\{x,y,z\}$ and grid spacing $h_d=L_d/M_d$. We use a separable window
\begin{equation}
W(\bm{r})=W_{\mathrm{pswf}}(x)\,W_{\mathrm{pswf}}(y)\,W_{\mathrm{pswf}}(z),
\end{equation}
with the one-dimensional kernel
\begin{equation}
W_{\mathrm{pswf}}(x):=\begin{cases}
\psi_0^{c}(x/\omega), \quad &|x|\leq \omega,\\[1.2mm]
0,\quad& \text{otherwise},
\end{cases}
\end{equation}
where $\omega=Ph/2$ is the half-width of the compact support in each direction, $h$ denotes the corresponding grid spacing, and $P\in\mathbb{Z}^{+}$ is the number of grid points coupled to each particle per dimension. The Fourier transform of $W_{\mathrm{pswf}}$ is available in closed form,
\begin{equation}
\widehat{W}_{\mathrm{pswf}}(k)=\omega\,\lambda_0\,\psi_{0}^{c}\!\left(\omega k/c\right),
\end{equation}
and is band-limited to $k\in[-c/\omega,c/\omega]$. Here, $\lambda_0$ is the eigenvalue associated with $\psi_0^{c}$ for the integral operator in \eqref{eq::formula}. Although $\psi_{0}^{c}$ has no elementary closed-form expression, it can be evaluated efficiently using piecewise polynomial approximations, so the cost is comparable to that of standard window functions.

The KB and ES windows can be interpreted as asymptotic approximations of PSWFs in appropriate parameter regimes~\cite{jiang2025dual,barnett_parallel_2019}. For a fixed target accuracy, KB and ES typically achieve the same error with smaller support than Gaussian and B-spline windows~\cite{shamshirgar2021fast}. PSWF windows often allow a further reduction in support at matched accuracy, thereby decreasing the number of grid points touched per particle during spreading~\cite{barnett_parallel_2019}. This advantage stems from the optimal frequency-concentration property of PSWFs under simultaneous compact-support and bandlimiting constraints~\cite{slepian1983sirev}, which improves aliasing-error control. A systematic comparison of window functions is left for future work.

\subsection{Implementation details}\label{subsec::implement}
We describe the implementation strategy for ESP, focusing on the spreading step and the efficient evaluation of PSWF-based kernels. For clarity, we first present the one-dimensional case and rescale the grid spacing to $1$, so that $W_{\mathrm{pswf}}(x)$ has support $[-P/2,P/2]$.

For a particle at position $x$, let $\bar{x}$ be the nearest grid point when $P$ is odd, or the midpoint between the two nearest grid points when $P$ is even. We then select the $P$ closest grid points
\begin{equation}
x_p:=\bar{x}+p,\qquad 
p\in\left\{-\left\lfloor\frac{P}{2}\right\rfloor,\,-\left\lfloor\frac{P}{2}\right\rfloor+1,\,\ldots,\,\left\lceil\frac{P}{2}\right\rceil-1\right\}.
\end{equation}
The weight assigned to grid point $x_p$ is
\begin{equation}
W_{p}(x):=W_{\mathrm{pswf}}\!\left(x-x_p\right)
= W_{\mathrm{pswf}}\!\left(x-(\bar{x}+p)\right),
\end{equation}
which depends only on the fractional offset $x-\bar{x}\in[-1/2,1/2]$.

Given a tolerance $\varepsilon$, we approximate $W_{\mathrm{pswf}}(x)$ by a piecewise polynomial $W_{\mathrm{poly}}(x)$ (e.g., via Chebyshev or Lagrange interpolation) on the $P$ unit-length subintervals of its support,
\begin{equation}
x\in\left[-\frac{P}{2}+ (n-1),\;-\frac{P}{2}+n\right],\qquad n=1,\ldots,P,
\end{equation}
so that
\begin{equation}
\|W_{\mathrm{poly}}-W_{\mathrm{pswf}}\|_{L^2}\le \varepsilon
\end{equation}
on each subinterval. In practice, we precompute and store the polynomial coefficients for each subinterval and reuse them for all particles. The resulting polynomials are evaluated efficiently using Horner's rule.

PSWFs also appear in the Fourier-space kernel (e.g., in $\widehat{\mathcal{F}}^{c}(\bm{k})$ and its derivative), and we use the same strategy of polynomial approximation when assembling the diagonal scaling for $\bm{P}_{\mathrm{ins},\mathcal{F}}$. For the short-range term $\bm{P}_{\mathrm{ins},\mathcal{N}}$ (Eq.~\eqref{eq::PinsN} truncated at $r_c$), we likewise approximate $\mathcal{F}_{\mathcal{N}}(r)$ to avoid repeated special-function evaluations.

To reduce the short-range cost further, we combine kernel precomputation with tabulation. Specifically, we introduce an inner cutoff $r_{\mathrm{in}}<r_c$. For $0<r\le r_{\mathrm{in}}$, we approximate $\mathcal{F}_{\mathcal{N}}(r)$ using a Taylor expansion. For $r_{\mathrm{in}}<r\le r_c$, we use a bitmask-based table lookup~\cite{wolff1999tabulated} with cubic-spline interpolation between table entries. This avoids expensive special-function calls while maintaining the prescribed accuracy.

Finally, the extension to three dimensions is straightforward because the window is separable: the one-dimensional weights are applied independently in each direction, yielding a $P^3$ stencil per particle for spreading and interpolation.

\subsection{Calculation of the local pressure tensor}\label{subsec::local}
Local pressure profiles (per particle or per group) are often needed to study inhomogeneous systems such as fluid interfaces, membranes, micelles, and self-assembled structures~\cite{ollila20093d,torres2015examining}. Extracting local contributions from the Fourier-space pressure (e.g., \eqref{eq::Far}) is nontrivial because the Fourier term is global rather than pairwise. For PME, this issue has been addressed in Ref.~\cite{sega2016pressure}. Here, we describe an analogous construction for ESP.

Our derivation starts from the end of Step~2 in Algorithm~\ref{al::Fourier}, where $\widehat{\rho}_{\mathrm{grid}}(\bm{k})$ is available. Define the mode-wise field
\begin{equation}\label{eq::rhodiag}
\widehat{\bm{\Pi}}^{\,\mathrm{local}}(\bm{k})
:=\widehat{\bm{\mathcal{P}}}^{c}(\bm{k})\,\widehat{W}(\bm{k})^{-2}\,\widehat{\rho}_{\mathrm{grid}}(\bm{k}),
\end{equation}
where $W$ is chosen real and even, so that $\widehat{W}(\bm{k})\in\mathbb{R}$. We then apply a single 3D inverse FFT to obtain $\bm{\Pi}^{\,\mathrm{local}}(\bm{r})$ on real-space grid points from $\widehat{\bm{\Pi}}^{\,\mathrm{local}}(\bm{k})$.

Using \eqref{eq::rhodiag} and Plancherel's theorem (\ref{app:PlancherelTheorem}), the long-range pressure contribution associated with particle $i$ can be written as
\begin{equation}\label{eq::PinsFind}
\begin{split}
\bm{P}_{\mathrm{ins},\mathcal{F}}^{\,i}
&=\frac{q_i}{2\det(\bm{h})^{2}}\sum_{\bm{k}\neq \bm{0}}
\widehat{W}(\bm{k})\,e^{i\bm{k}\cdot\bm{r}_i}\,
\widehat{\bm{\Pi}}^{\,\mathrm{local}}(\bm{k})\\
&=\frac{q_i}{2}\int_{\Omega}\bm{\Pi}^{\,\mathrm{local}}(\bm{r})\,W_{*}(\bm{r}_i-\bm{r})\,d\bm{r},
\end{split}
\end{equation}
where $W_{*}$ denotes the periodized window. Discretizing the integral with the trapezoidal rule and evaluating it at particle positions yields a standard gathering step. Since $W$ is smooth and compactly supported, this procedure converges spectrally as the grid is refined.

Compared with the four-step global evaluation in Section~\ref{subsec::FFTSpectral}, computing local contributions adds two operations: an inverse FFT and a gather. In practice, this overhead is modest because local pressure profiles are typically accumulated for a subset of configurations (e.g., during sampling or ensemble averaging). Moreover, when the same window is used for spreading and gathering, the stencil indices and window values can be cached and reused to further reduce cost.

\subsection{Complexity analysis}\label{subsec::complexity}
We summarize the computational complexity of the ESP pressure evaluation, which consists of a short-range (real-space) contribution and a long-range (Fourier-space) contribution.

\noindent \textbf{Short-range cost.}
The short-range pressure $\bm{P}_{\mathrm{ins},\mathcal{N}}$ is evaluated by direct truncation at the cutoff radius $r_c$ using a cell-list (Verlet list) approach~\cite{verlet1967computer}. Its cost scales as
\[
O(n_s N),
\]
where
\[
n_s=\frac{4\pi}{3}r_c^3\rho_r
\]
is the average number of particles within distance $r_c$ at number density $\rho_r$.

\noindent \textbf{Long-range cost.}
Let $I_{\mathcal{F}}$ denote the number of reciprocal-space modes used in the long-range evaluation. In ESP, the reciprocal-space kernel is effectively supported on $|\bm{k}|\le K_{\max}$ with $K_{\max}=c/r_c$ (cf.\ \eqref{eq::Far} and the compact support of $\psi_0^c$ and $\psi_0^{c}{}^{\prime}$ on $[-1,1]$).
Because the reciprocal-space sums are evaluated on an FFT mesh, the set of represented wavevectors is naturally bounded by a \emph{rectangular} (box-shaped) truncation aligned with the reciprocal lattice directions rather than a spherical cutoff. Geometrically, the number of included modes can therefore be estimated as
\begin{equation}\label{eq::IF_geom}
I_{\mathcal{F}}
\;\approx\;
\frac{\mathrm{vol}\bigl(\{\bm{k}:\ |k_i|\le K_{i,\max},\ i=1,2,3\}\bigr)}{\mathrm{vol}(\mathcal{B}^*)}
\;=\;
\frac{8\,K_{1,\max}K_{2,\max}K_{3,\max}}{|\det[\bm{b}_1\,\bm{b}_2\,\bm{b}_3]|},
\end{equation}
where $\mathcal{B}^*$ is the reciprocal lattice cell spanned by $\bm{b}_1,\bm{b}_2,\bm{b}_3$ (defined in \eqref{eq::vector}). Using the identity
\[
|\det[\bm{b}_1\,\bm{b}_2\,\bm{b}_3]|=\frac{(2\pi)^3}{\det(\bm{h})}=\frac{(2\pi)^3}{V},
\]
we obtain the explicit geometric estimate
\begin{equation}\label{eq::IF_const}
I_{\mathcal{F}}
\;\approx\;
\frac{V}{(2\pi)^3}\,8\,K_{1,\max}K_{2,\max}K_{3,\max}
\;=\;
\frac{V}{\pi^3}\,K_{1,\max}K_{2,\max}K_{3,\max}.
\end{equation}
In the common isotropic setting $K_{1,\max}=K_{2,\max}=K_{3,\max}=K_{\max}$ (e.g., when using the same cutoff parameter in each direction), this reduces to
\begin{equation}\label{eq::IF_iso}
I_{\mathcal{F}}
\;\approx\;
\frac{V}{\pi^3}\,K_{\max}^3
\;=\;
\frac{N}{\pi^3\rho_r}\left(\frac{c}{r_c}\right)^3.
\end{equation}
Equivalently, in terms of $n_s=\frac{4\pi}{3}r_c^3\rho_r$,
\begin{equation}\label{eq::IF_ns}
I_{\mathcal{F}}
\;\approx\;
\frac{4}{3\pi^2}\,\frac{c^3}{n_s}\,N.
\end{equation}
In practice, the FFT mesh corresponds to a rectangular index set in reciprocal space; one may therefore view \eqref{eq::IF_const}--\eqref{eq::IF_ns} as a geometric baseline, with the actual constant depending mildly on the chosen mesh shape and any additional oversampling.

Using \eqref{eq::c}, we have
\[
I_{\mathcal{F}}
=O\!\left(\log^3(\varepsilon^{-1})\,r_c^{-3}\rho_r^{-1}\,N\right)
=O\!\left(\log^3(\varepsilon^{-1})\,n_s^{-1}\,N\right).
\]
In Algorithm~\ref{al::Fourier}, spreading couples each particle to $P^3$ grid points, so its cost is $O(P^3N)$. The diagonal scaling and Fourier-space reduction cost $O(I_{\mathcal{F}})$, while the forward 3D FFT costs $O(I_{\mathcal{F}}\log I_{\mathcal{F}})$.

\noindent \textbf{Total complexity.}
Combining these costs, the per-step complexity of ESP is
\begin{equation}\label{eq::complexity}
C_{\mathrm{ESP}}
=O\!\left(n_sN + P^3N + I_{\mathcal{F}}\log I_{\mathcal{F}}\right)
=O\!\left(n_sN + P^3N + \log^3(\varepsilon^{-1})\,n_s^{-1}N\,\log N\right).
\end{equation}
In practice, one may either keep $r_c=O(1)$ (so that $\rho_r$ is fixed and hence $n_s=O(1)$), in which case fixing $P$ by the accuracy target recovers the usual $O(N\log N)$ scaling, or tune $r_c$ to balance the real- and reciprocal-space costs at a prescribed tolerance $\varepsilon$. In particular, since the real-space work scales like $n_sN$ while the reciprocal-space work contributes a term of the form $\log^3(\varepsilon^{-1})\,n_s^{-1}N$ (through the required Fourier resolution), choosing
\begin{equation}
n_s = O\!\left(\log^{3/2}(\varepsilon^{-1})\right)
\end{equation}
balances these contributions and yields an overall $O(N\log N)$ complexity with a leading prefactor scaling like $\log^{3/2}(\varepsilon^{-1})$. The optimal balance is implementation-dependent, however, as it also reflects the relative constants in the direct short-range evaluation and the FFT/pencil-decomposition communication.

For the local-pressure procedure in Section~\ref{subsec::local}, we add one inverse FFT and one gathering step, with costs $O(I_{\mathcal{F}}\log I_{\mathcal{F}})$ and $O(P^3N)$, respectively. The overall complexity therefore remains $O(N\log N)$.

\section{Numerical results}\label{sec::NumericalExample}
In this section, we present numerical tests to assess the accuracy and efficiency of ESP-based molecular dynamics in the NPT ensemble. The test systems include SPC/E bulk water, LiTFSI ionic liquids, and a bovine cytochrome bc$_1$ transmembrane complex. Our implementations are built on the open-source packages LAMMPS (19 Nov 2024) and GROMACS (2025.1). All calculations were performed on the Flatiron Institute ``Rusty'' cluster supported by the Scientific Computing Core. Each node has two AMD EPYC 9474F 48-core CPUs and 1.5~TB of memory. The implementations are freely available (see Refs.~\cite{lammpsespcode,gromacsespcode}), and the accompanying scripts reproduce all numerical results; timing measurements are expected to agree up to machine- and system-dependent variability.

\subsection{Parameter selection and accuracy on bulk SPC/E water}
We first evaluate the accuracy of ESP for NPT simulations of SPC/E water~\cite{Berendsen1987JPCSPCE} using LAMMPS. The system contains $21{,}624$ atoms in a cubic box with side length $3~\mathrm{nm}$. Each run starts with $200~\mathrm{ps}$ of NPT equilibration using anisotropic pressure coupling at $300~\mathrm{K}$ and a reference pressure of $10^{-3}~\mathrm{kbar}$, with time step $\Delta t=2~\mathrm{fs}$. Temperature and pressure are controlled using the Martyna--Tobias--Klein (MTK) algorithm~\cite{martyna1994constant} with a relaxation time of $0.1~\mathrm{ps}$. We then run $200~\mathrm{ps}$ of production and save configurations every $2~\mathrm{ps}$ (100 snapshots) for analysis. Pressure-tensor errors are computed using relative $L^2$ norms over these saved configurations. We use a real-space cutoff of $r_c=0.9~\mathrm{nm}$ for all methods.

We quantify errors in the pressure tensor by the relative $L^2$ norms of the diagonal and off-diagonal components. Define
\[
e_{\alpha\beta}:=(\bm{P}_{\mathrm{ins}}^{\mathrm{ref}})_{\alpha\beta}-(\bm{P}_{\mathrm{ins}}^{\mathrm{esp}})_{\alpha\beta}.
\]
Then
\begin{equation}\label{eq::ref_pres}
\mathcal{E}_{\mathrm{diag}}
:=\frac{\left(\sum_{\alpha=1}^{3} e_{\alpha\alpha}^2\right)^{1/2}}
{\left(\sum_{\alpha=1}^{3} (\bm{P}_{\mathrm{ins}}^{\mathrm{ref}})_{\alpha\alpha}^2\right)^{1/2}},
\qquad
\mathcal{E}_{\mathrm{off\mbox{-}diag}}
:=\frac{\left(\sum_{1\le \alpha<\beta\le 3} e_{\alpha\beta}^2\right)^{1/2}}
{\left(\sum_{1\le \alpha<\beta\le 3} (\bm{P}_{\mathrm{ins}}^{\mathrm{ref}})_{\alpha\beta}^2\right)^{1/2}}.
\end{equation}
We compute the reference pressure tensor $\bm{P}_{\mathrm{ins}}^{\mathrm{ref}}$ using PPPM~\cite{Hockney1988Computer} with tolerance $10^{-12}$.

For a target tolerance $\Delta$, we choose the PSWF splitting and spreading parameters $c_{\mathrm{split}}$ and $c_{\mathrm{spread}}$ such that
\[
\psi_0^{c_{\mathrm{split}}}(1)=\psi_0^{c_{\mathrm{spread}}}(1)=\Delta.
\]
In the high-precision regime (small $\Delta$), this implies the scaling $c\approx \log(1/\Delta)$. We approximate $\psi_0^{c_{\mathrm{split}}}$ and $\psi_0^{c_{\mathrm{spread}}}$ using 18th-order polynomial interpolation; the interpolation error is sufficiently small that it does not affect the reported accuracy. The spreading order $P$ is taken identical in all dimensions. Although the box size fluctuates slightly over the trajectory, these parameters remain adequate and are kept fixed throughout the simulation.

\begin{figure}[!htbp]
\centering
\includegraphics[width=0.98\linewidth]{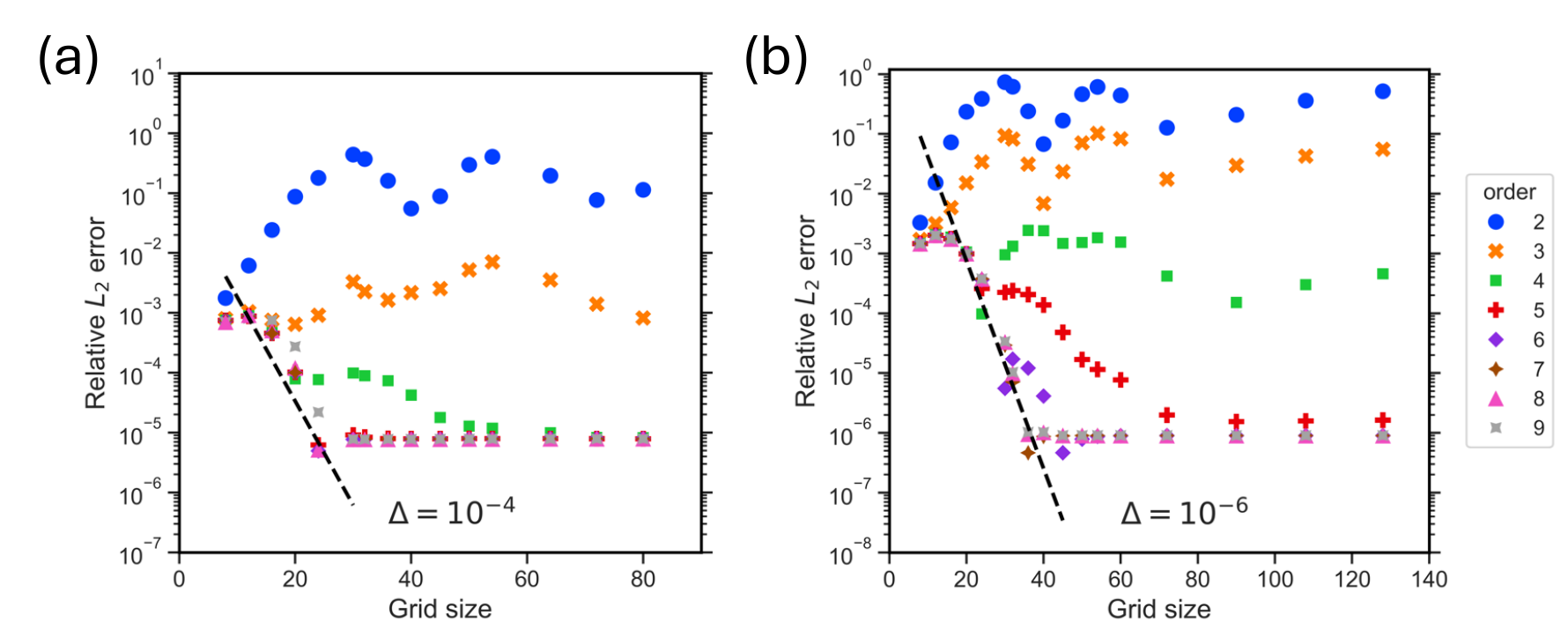}
\caption{\sf Relative $L_2$ error of the diagonal components of $\bm{P}_{\mathrm{ins}}$, computed over 100 equilibrium configurations, as a function of the per-dimension grid size $I_d$ (total grid points $I_{\mathcal{F}}=I_d^3$). Error tolerances are {\bf (a)} $\Delta=10^{-4}$ and {\bf (b)} $\Delta=10^{-6}$. Results are shown for several spreading orders $P$. The dashed line shows the fitted convergence rate, $O(e^{-I_d/2.5})$.}
\label{fig:diagonalaccuracy}
\end{figure}

\begin{figure}[!htbp]
\centering
\includegraphics[width=0.98\linewidth]{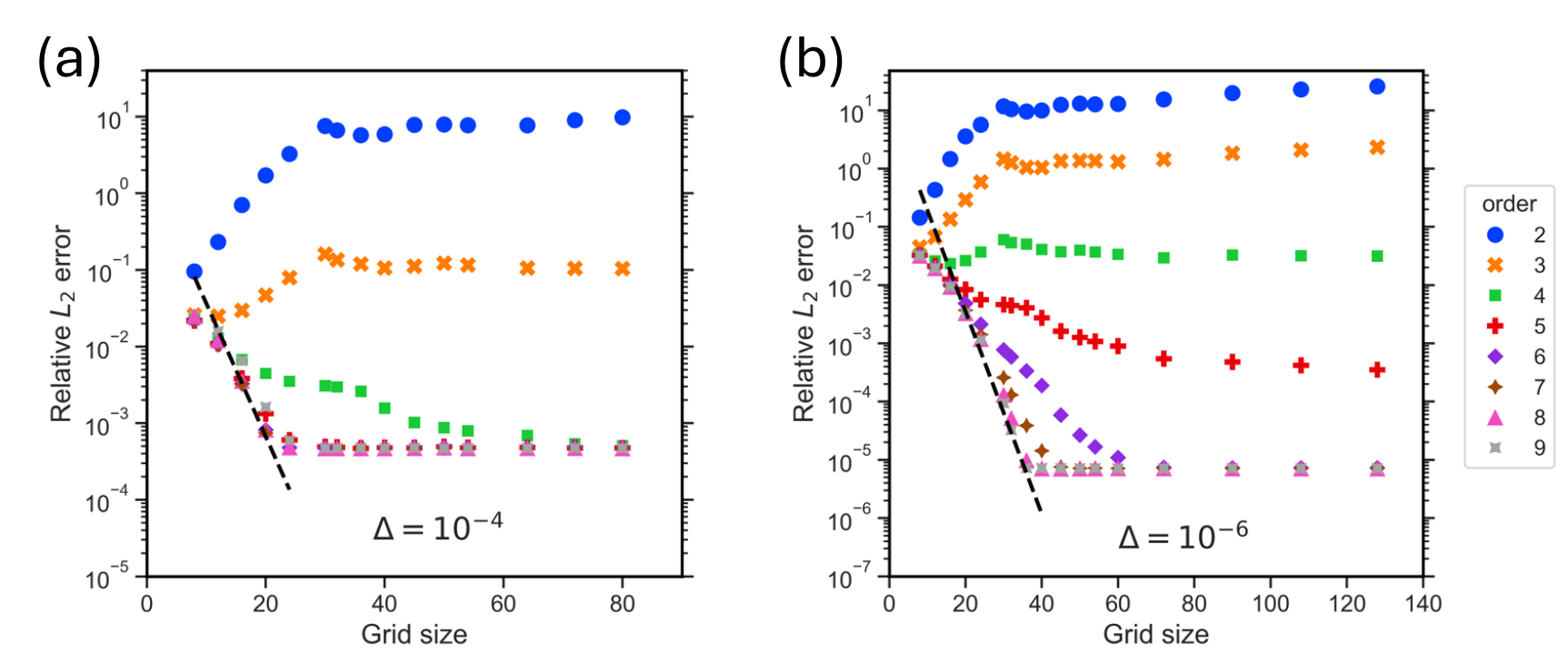}
\caption{\sf Relative $L_2$ error of the off-diagonal components of $\bm{P}_{\mathrm{ins}}$, computed over 100 equilibrium configurations, as a function of the per-dimension grid size $I_d$. Error tolerances are {\bf (a)} $\Delta=10^{-4}$ and {\bf (b)} $\Delta=10^{-6}$. Results are shown for several spreading orders $P$. The dashed line shows the fitted convergence rate, $O(e^{-I_d/2.5})$.}
\label{fig:offdiagonalaccuracy}
\end{figure}

Figures~\ref{fig:diagonalaccuracy} and \ref{fig:offdiagonalaccuracy} plot $\mathcal{E}_{\mathrm{diag}}$ and $\mathcal{E}_{\mathrm{off\mbox{-}diag}}$, respectively, as functions of the per-dimension grid size $I_d$ (so the total number of grid points is $I_{\mathcal{F}}=I_d^3$). Results are shown for several spreading orders $P$ and two target tolerances, $\Delta=10^{-4}$ and $10^{-6}$. The results exhibit spectral convergence as the grid is refined. They also indicate that, for a target tolerance $\Delta=10^{-D}$, one typically needs $P\ge D+1$ to avoid excessive Fourier upsampling: when $P\le D$, achieving the target accuracy requires a substantial increase in the number of Fourier modes; when $P=D+1$, only mild upsampling is needed; and when $P>D+1$, upsampling is typically unnecessary. Thus, good performance requires balancing $P$ against the Fourier grid size. Unlike the standard PPPM implementation in LAMMPS, which supports only $P\le 7$, ESP supports arbitrary $P$ and can, in principle, reach arbitrarily high accuracy.

\begin{figure}[!htbp]
\centering
\includegraphics[width=0.98\linewidth]{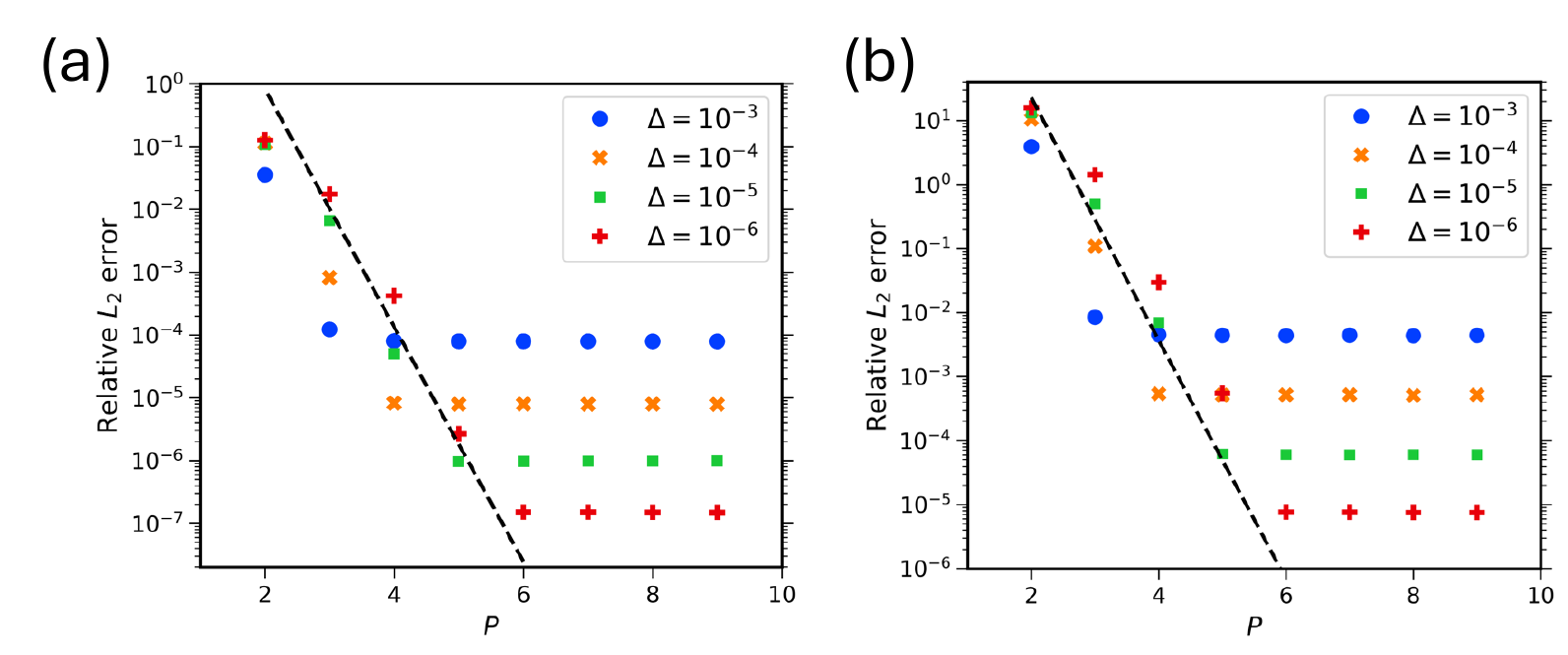}
\caption{\sf Relative $L_2$ error of the {\bf (a)} diagonal and {\bf (b)} off-diagonal components of $\bm{P}_{\mathrm{ins}}$, computed over 100 equilibrium configurations, as a function of the spreading order $P$. Results are shown for different tolerances $\Delta$. The dashed line shows the fitted convergence rate, $O(\mathrm{erfc}(2.05\sqrt{P}))$.}
\label{fig:accuracyspreadingorder}
\end{figure}

To further assess how the spreading order affects accuracy, we evaluate the pressure-tensor error for varying $P$ and target tolerances $\Delta\in[10^{-3},10^{-6}]$, while fixing the Fourier grid to $I_{\mathcal{F}}=80^3$. With this choice, Fourier truncation errors are negligible and the observed error is dominated by spreading. Figures~\ref{fig:accuracyspreadingorder}(a)--(b) show $\mathcal{E}_{\mathrm{diag}}$ and $\mathcal{E}_{\mathrm{off\mbox{-}diag}}$ as functions of $P$. The errors decay rapidly and are well fit by $O(\mathrm{erfc}(C_1\sqrt{P}))$, with fitted constant $C_1=2.05$ for both diagonal and off-diagonal components. Compared with the decay rate reported for Gaussian windows~\cite{gao_fast_2025}, this corresponds to an approximately $1.7\times$ faster decay in $P$. These results suggest that relatively small $P$ is sufficient for most MD applications, where the target tolerance is typically no tighter than $10^{-5}$.

\begin{figure}[!htbp]
\centering
\includegraphics[width=0.98\linewidth]{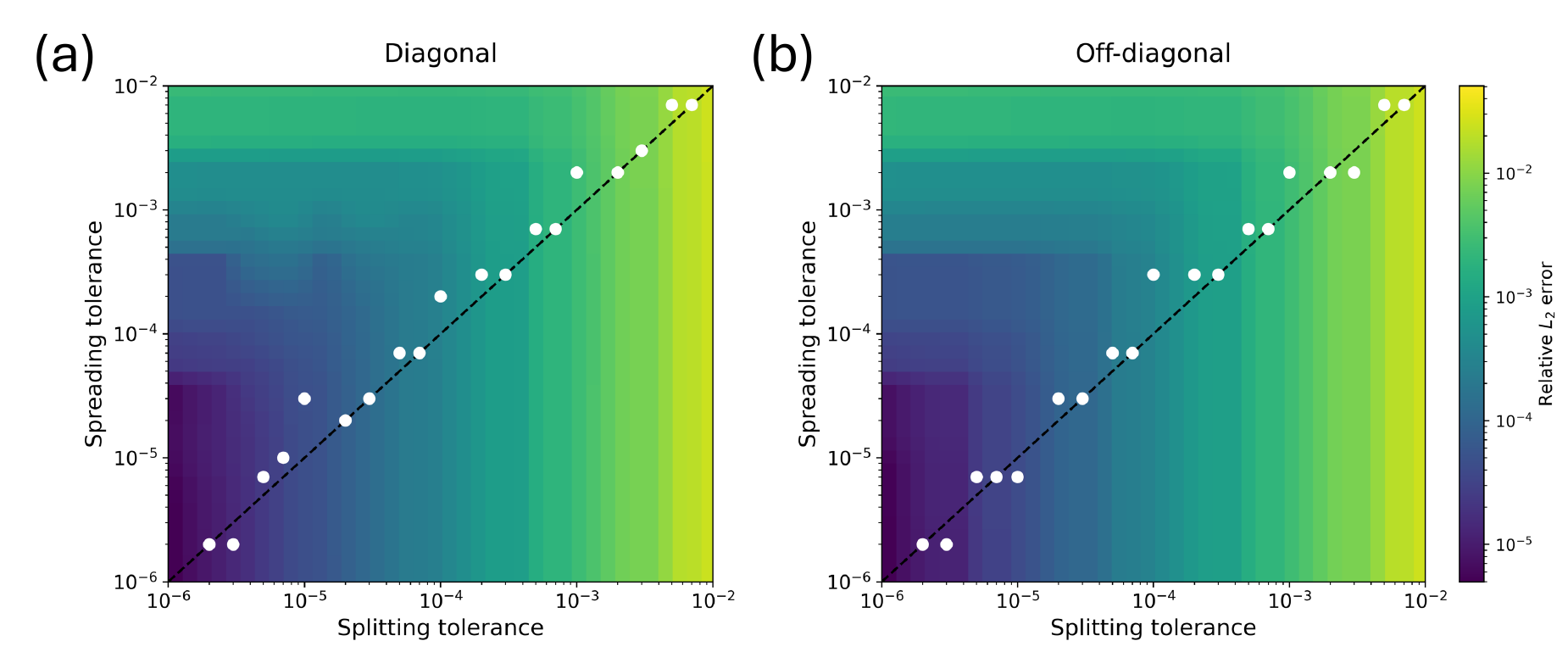}
\caption{\sf Heatmaps of the relative error in the {\bf (a)} diagonal and {\bf (b)} off-diagonal components of the pressure tensor. The $x$- and $y$-axes denote the splitting tolerance $\Delta_{\mathrm{split}}$ and spreading tolerance $\Delta_{\mathrm{spread}}$, respectively. White circles mark the ``optimal'' pairs, defined as the largest $\Delta_{\mathrm{spread}}$ that achieves an error within $1.01\times$ the minimum error attained at that $\Delta_{\mathrm{split}}$.}
\label{fig:accuracyHeatMap}
\end{figure}

In the previous tests, we set the splitting and spreading tolerances to be equal. To assess this choice, we study how the pressure-tensor error depends on $c_{\mathrm{split}}$ and $c_{\mathrm{spread}}$. Figures~\ref{fig:accuracyHeatMap}(a)--(b) show heatmaps of the relative $L_2$ errors for the diagonal and off-diagonal components, respectively. The $x$-axis is the splitting tolerance $\Delta_{\mathrm{split}}$ and the $y$-axis is the spreading tolerance $\Delta_{\mathrm{spread}}$, where $c_{\mathrm{split}}$ and $c_{\mathrm{spread}}$ are chosen such that
$\psi_{0}^{c_{\mathrm{split}}}(1)=\Delta_{\mathrm{split}}$ and $\psi_{0}^{c_{\mathrm{spread}}}(1)=\Delta_{\mathrm{spread}}$.
In these tests, we fix $P=7$ and $I_{\mathcal{F}}=80^3$ to isolate the effect of $(c_{\mathrm{split}},c_{\mathrm{spread}})$ on accuracy. To build the heatmaps, we first evaluate errors on a $20\times 20$ grid, sampling both tolerances at logarithmically spaced values from $10^{-2}$ to $10^{-6}$. We then apply bilinear interpolation in the $(\log \Delta_{\mathrm{split}},\log \Delta_{\mathrm{spread}})$ plane to obtain a smooth $40\times 40$ grid for visualization. For each $\Delta_{\mathrm{split}}$, we define an ``optimal'' choice of $\Delta_{\mathrm{spread}}$ as the largest value whose error is within $1.01\times$ the minimum error achieved at that $\Delta_{\mathrm{split}}$ (a 1\% buffer is included to reduce sensitivity to rounding). The resulting optimal pairs cluster near the diagonal, which supports using $c_{\mathrm{split}}=c_{\mathrm{spread}}$ as a simple and near-optimal choice over a wide range of target accuracies. Based on these results, we set $c_{\mathrm{split}}=c_{\mathrm{spread}}$ in all subsequent simulations and do not mention this choice again.

\begin{figure}[!htbp]
\centering
\includegraphics[width=0.65\linewidth]{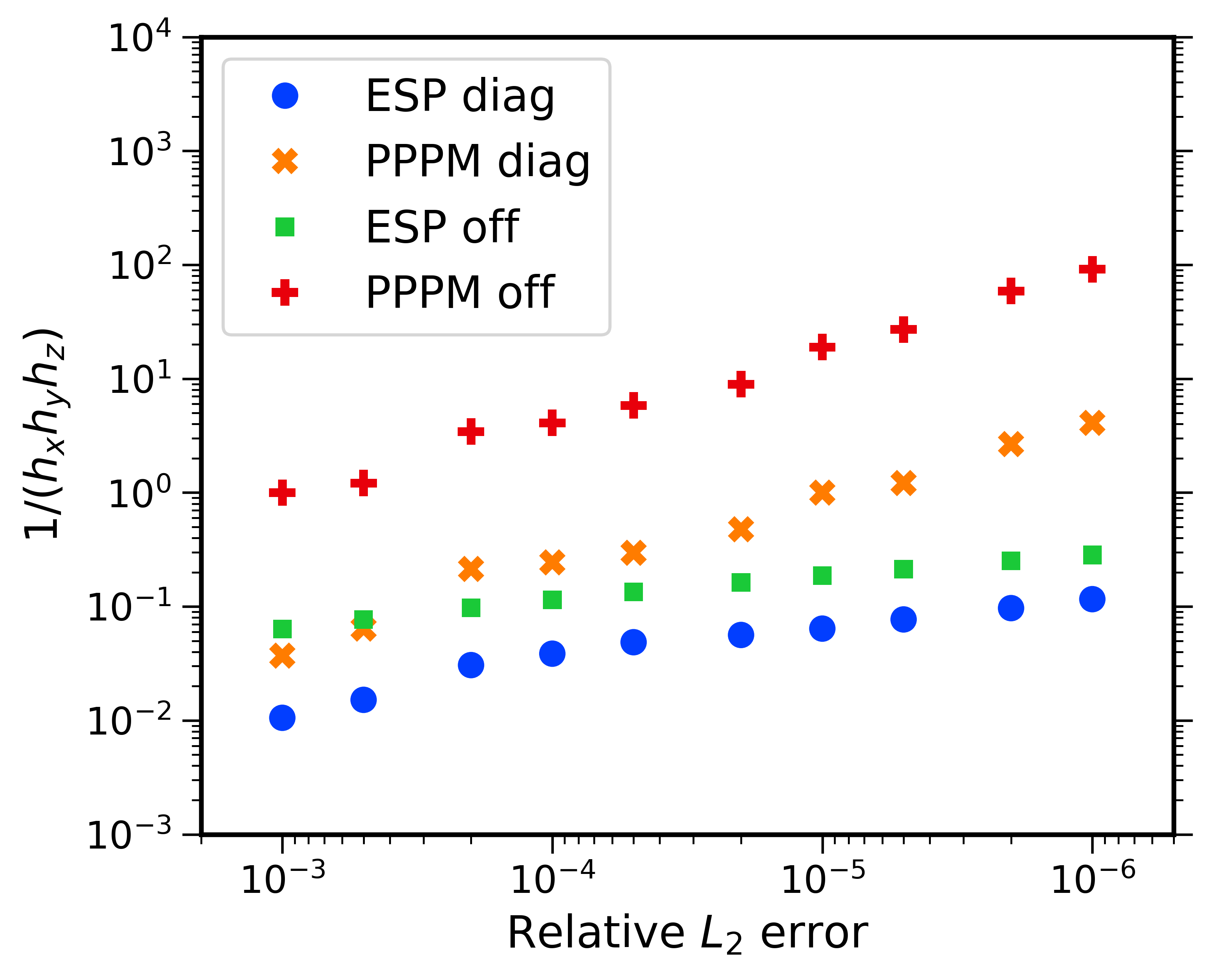}
\caption{\sf Relative $L_2$ error of the instantaneous pressure tensor versus the smallest inverse mesh volume $1/(h_xh_yh_z)$ required to reach the specified error level. For ESP, circles and squares denote the errors of the diagonal and off-diagonal components, respectively. For PPPM, crosses and plus signs denote the corresponding diagonal and off-diagonal errors.}
\label{fig:GridAcc}
\end{figure}

Next, we compare ESP and PPPM by plotting the relative $L_2$ errors of the diagonal and off-diagonal pressure-tensor components against the inverse mesh volume $1/(h_xh_yh_z)$, which is directly proportional to the number of Fourier modes. For PPPM, we use the native LAMMPS implementation with the recommended settings: spreading order $P=5$, and $h_x$, $h_y$, and $h_z$ are selected automatically from the target tolerance and the corresponding error-bound estimate~\cite{deserno1998meshII}. For ESP, we set $\Delta=\Delta_{\mathrm{split}}=\Delta_{\mathrm{spread}}$ and vary $\Delta$ to match the target error levels. We choose $P=\lceil-\log_{10}(\Delta)\rceil+1$, consistent with the trends in Figs.~\ref{fig:diagonalaccuracy}--\ref{fig:offdiagonalaccuracy}. For example, to reach relative $L_2$ errors of $10^{-3}$, $10^{-4}$, and $10^{-5}$, we use $\Delta=2\times10^{-2}$, $10^{-3}$, and $2\times10^{-4}$ for the diagonal component, and $\Delta=2\times10^{-4}$, $2\times10^{-5}$, and $2\times10^{-6}$ for the off-diagonal component. Figure~\ref{fig:GridAcc} shows that ESP attains the same error at substantially smaller $1/(h_xh_yh_z)$ than PPPM, indicating that ESP can use a significantly coarser mesh. The mesh-point reduction increases as the target error decreases: over the relative $L_2$ error range $10^{-3}$ to $10^{-6}$, the reduction is $3.5$--$35.2\times$ for the diagonal component and $15.7$--$322.1\times$ for the off-diagonal component. Because $P$ increases as $\Delta$ decreases, the wall-time speedup for pressure evaluation may be smaller than the mesh-point reduction alone suggests. However, FFT all-to-all communication scales with the number of Fourier modes, so ESP is expected to substantially reduce communication cost in large-scale NPT simulations; we verify this in the next section.

\subsection{Performance and strong-scaling benchmarks}
We evaluate the CPU performance of ESP for NPT simulations of SPC/E bulk water~\cite{Berendsen1987JPCSPCE}. We first compare against the PME implementation in GROMACS using its recommended settings (spreading order $P=5$, cutoff $r_c=0.9~\mathrm{nm}$, and Fourier spacing $0.12~\mathrm{nm}$). With these settings, the measured relative errors are $4\times10^{-4}$ for forces and $2\times10^{-5}$ (diagonal) and $2\times10^{-3}$ (off-diagonal) for the pressure-tensor components. For ESP, we choose $\Delta=4\times10^{-4}$, $P=5$, $r_c=0.9~\mathrm{nm}$, and Fourier spacing $0.26~\mathrm{nm}$, which matches the same accuracy level for both forces and the pressure tensor. All timings are averaged over 1{,}000 steps.

\begin{figure}[!htbp]
\centering
\includegraphics[width=0.98\linewidth]{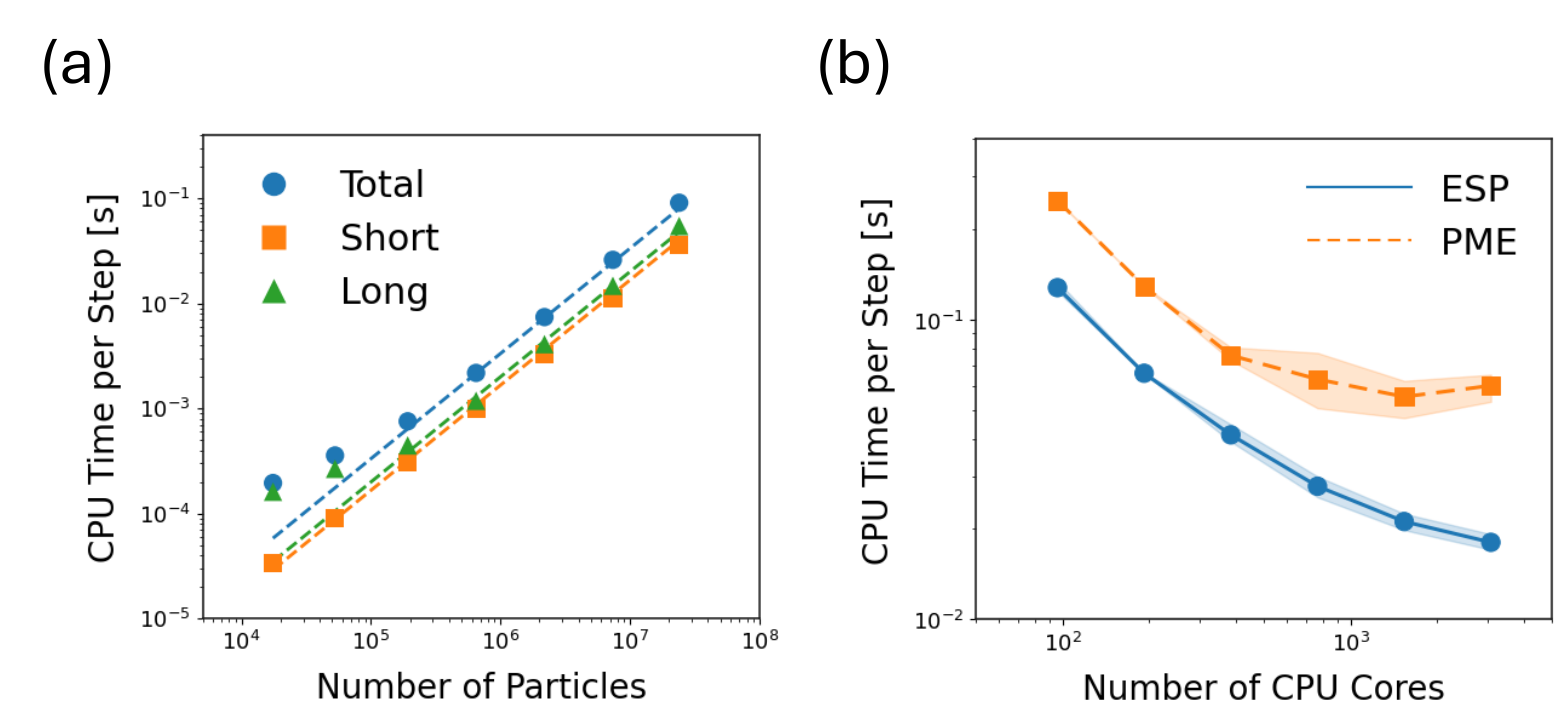}
\caption{\sf Wall-clock time per step for Coulomb calculations (forces and pressure tensor) versus {\bf (a)} particle number and {\bf (b)} number of CPU cores. In {\bf (a)}, we report the short-range, long-range, and total costs for ESP with $\Delta=4\times10^{-4}$, with the number of CPU cores fixed at 96; dashed lines show linear fits. In {\bf (b)}, the particle number is fixed at $24{,}000{,}000$, and we compare the total runtime of ESP and PME (GROMACS). Shaded regions indicate confidence intervals estimated from five repeated runs.}
\label{fig:performance}
\end{figure}

Figure~\ref{fig:performance}(a) reports the wall-clock time per step for Coulomb calculations as the particle number increases from $N=17{,}496$ to $24{,}000{,}000$ with 96 CPU cores. Because the NPT update requires both forces and the instantaneous pressure tensor, we evaluate them together (see \ref{app::force}) and report the combined cost. The short-range, long-range, and total costs of ESP grow approximately linearly with $N$. This behavior is expected in the present regime: although the spectral solver has $O(N\log N)$ complexity, the FFT cost is not dominant at these resolutions because charge spreading and short-range interactions account for most of the runtime and scale as $O(N)$. Figure~\ref{fig:performance}(b) compares strong scaling of the total Coulomb time for ESP and PME on a system with $N=24{,}000{,}000$ particles. ESP increasingly outperforms PME as the core count grows, and the gap becomes more pronounced at large scale. At 3{,}072 cores, ESP is approximately $3\times$ faster than PME. Overall, these results demonstrate strong parallel scalability of ESP for large NPT simulations.

Next, we consider a higher-accuracy setting and compare against the PPPM implementation in LAMMPS. For ESP, we use $\Delta=2\times10^{-5}$, $P=6$, $r_c=0.9~\mathrm{nm}$, and Fourier spacing $0.20~\mathrm{nm}$. The measured relative errors are about $2\times10^{-5}$ for forces and $1\times10^{-6}$ (diagonal) and $1\times10^{-4}$ (off-diagonal) for the pressure components. For PPPM, we use the LAMMPS default spreading order $P=5$, $r_c=0.9~\mathrm{nm}$, and Fourier spacing $0.067~\mathrm{nm}$ to match the same accuracy level for both forces and the pressure tensor. With these settings, ESP reduces the total number of Fourier modes by approximately $27\times$. Timings are averaged over 1{,}000 steps.

\begin{figure}[!htbp]
\centering
\includegraphics[width=0.98\linewidth]{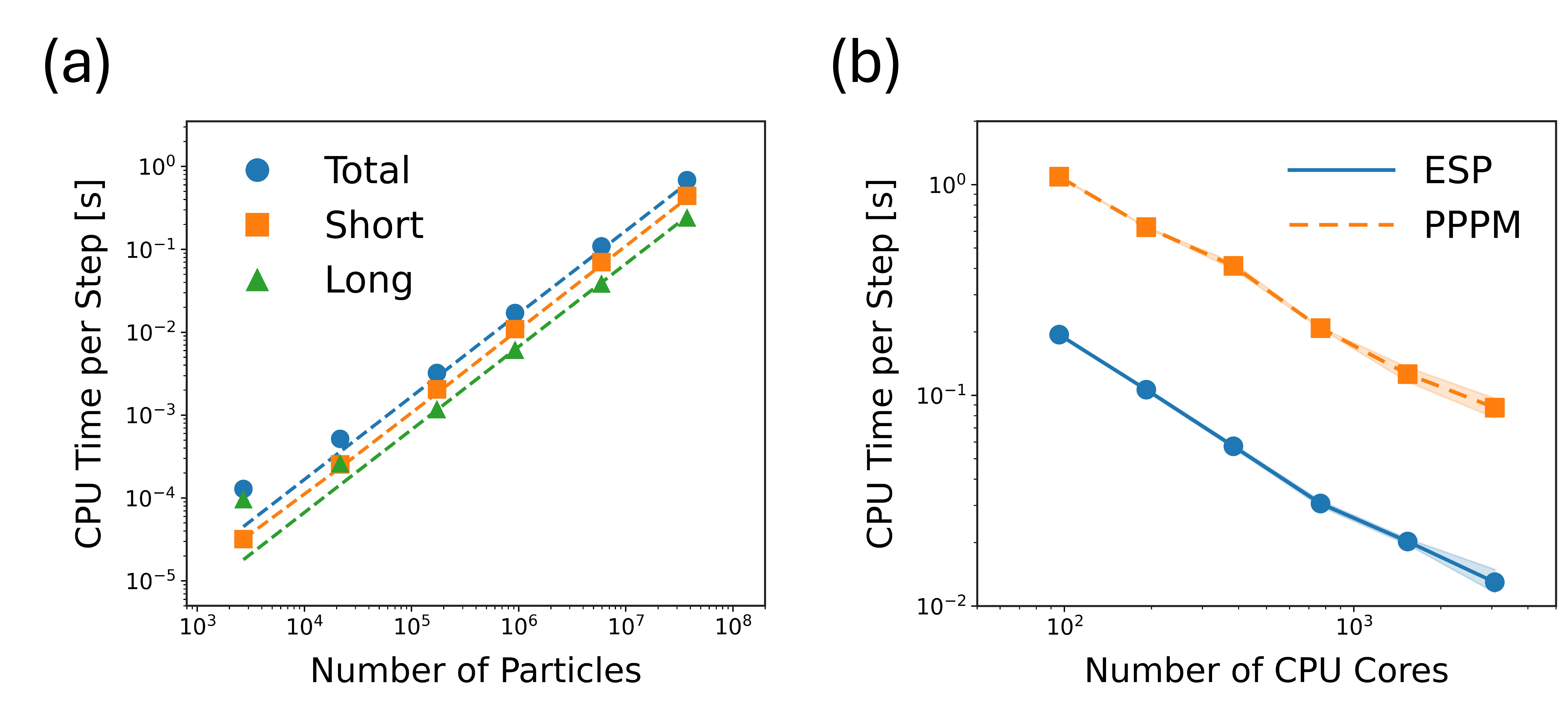}
\caption{\sf High-accuracy CPU performance test (ESP vs PPPM) for NPT simulations of SPC/E bulk water. Wall-clock time per step for Coulomb calculations (forces and pressure tensor) versus {\bf (a)} particle number and {\bf (b)} number of CPU cores. In {\bf (a)}, we report the short-range, long-range, and total costs for ESP with $\Delta=2\times10^{-5}$, with the number of CPU cores fixed at 96; dashed lines show linear fits. In {\bf (b)}, the particle number is fixed at $11{,}071{,}488$, and we compare the total runtime of ESP and PPPM (LAMMPS). Shaded regions indicate confidence intervals estimated from five repeated runs.}
\label{fig:performancehighacc}
\end{figure}

Figures~\ref{fig:performancehighacc}(a)--(b) report the Coulomb time per step versus particle number (from $N=2{,}703$ to $37{,}366{,}000$ at 96 cores) and the strong-scaling comparison on a system with $N=11{,}071{,}488$ particles, respectively. We again observe close-to-linear growth with $N$. In the strong-scaling test, ESP is about $5$--$8\times$ faster than PPPM, with a larger speedup than in the lower-accuracy case. While this accuracy level may not be required for routine MD, it can be valuable in pressure-sensitive settings such as high-pressure fluids, phase equilibrium and phase transitions, and local stress analysis near nanopores.

\subsection{Long-time simulations of LiTFSI ionic liquid}
\begin{figure}[!htbp]
\centering
\includegraphics[width=0.92\linewidth]{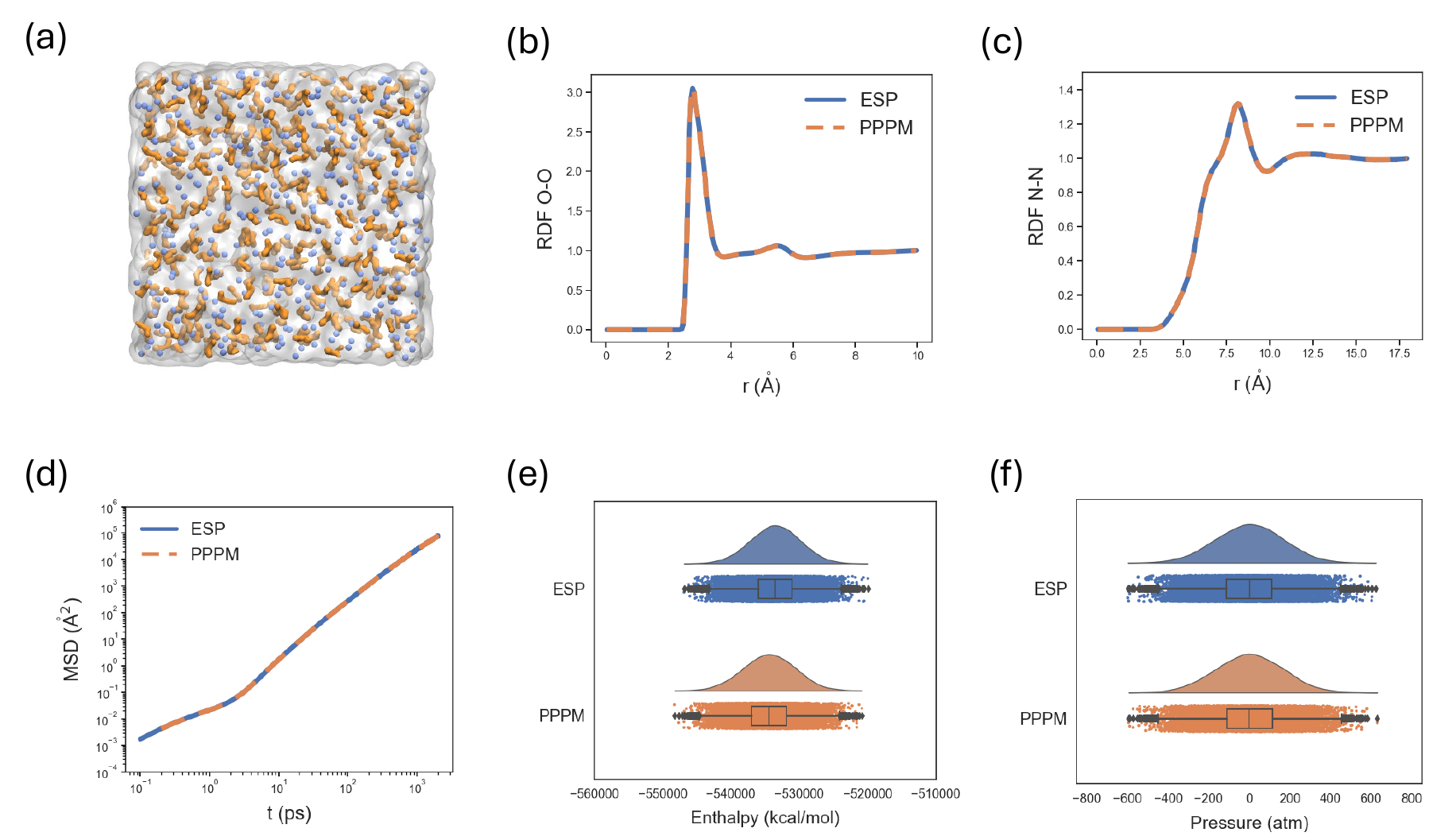}
\caption{\sf Comparison of PPPM and ESP for a concentrated aqueous LiTFSI electrolyte ($126{,}424$ atoms). Both methods target a relative error tolerance of $10^{-4}$. {\bf (a)} Snapshot of the system: anions (yellow), cations (blue), and water (white; rendered with QuickSurf). {\bf (b,c)} Radial distribution functions (RDFs) for nitrogen--nitrogen (N--N) pairs between anions and oxygen--oxygen (O--O) pairs between water molecules, respectively. {\bf (d)} Mean-squared displacement of nitrogen atoms. {\bf (e,f)} Raincloud plots of the total enthalpy and pressure, showing the probability density together with a boxplot and individual samples. PPPM and ESP yield statistically consistent distributions.}
\label{fig:ionicliquid}
\end{figure}

We benchmark ESP on a concentrated aqueous LiTFSI electrolyte at an ultrahigh concentration of $5~\mathrm{mol/L}$. The system contains $126{,}424$ atoms, including $2{,}560$ Li$^{+}$, $2{,}560$ TFSI$^{-}$, and $28{,}488$ H$_2$O molecules. Simulations are performed in a cubic box of side length $11.46~\mathrm{nm}$ with isotropic pressure coupling. At this concentration, the electrolyte is microscopically inhomogeneous and forms water-rich and anion-rich domains on $1$--$2~\mathrm{nm}$ length scales; both networks percolate through the simulation cell. Figure~\ref{fig:ionicliquid}(a) shows an MD snapshot illustrating this nano-heterogeneity. We use our LAMMPS implementation of ESP and compare against the native PPPM method~\cite{Hockney1988Computer}. We model LiTFSI using OPLS-AA parameters for Li$^{+}$, TIP3P for water~\cite{price2004modified}, and a molecular force field for TFSI$^{-}$~\cite{canongia2004molecular}. We equilibrate the system in the NPT ensemble at $298~\mathrm{K}$ and $1~\mathrm{bar}$ for $200~\mathrm{ns}$, and then run an additional $200~\mathrm{ns}$ of NPT production using the Martyna--Tobias--Klein (MTK) barostat~\cite{martyna1994constant} with either ESP or PPPM.

In all simulations, the cutoff for short-range Coulomb and Lennard--Jones interactions is set to $r_c=0.9~\mathrm{nm}$, and we use spreading order $P=5$ (LAMMPS default). For ESP, we set $\Delta=2\times10^{-4}$ for both splitting and spreading, with mesh spacing $0.24~\mathrm{nm}$. For PPPM, the mesh spacing is refined to $0.11~\mathrm{nm}$ to obtain a comparable relative error. To assess accuracy, we compute the radial distribution functions (RDFs) for N--N and O--O pairs and the mean-squared displacement (MSD) of nitrogen atoms. The RDFs characterize liquid structure, while the MSD captures translational motion over a wide range of time scales. The agreement in Fig.~\ref{fig:ionicliquid}(b)--(d) shows that ESP reproduces both structural and dynamical properties while using about $1/12$ as many Fourier modes.

We further assess statistical consistency by comparing the mean values and fluctuations of enthalpy and pressure. The time-averaged distributions are shown in Fig.~\ref{fig:ionicliquid}(e)--(f), and both methods produce indistinguishable statistics within sampling uncertainty. On 384 CPU cores, ESP achieves $124.58~\mathrm{ns/day}$, compared with $43.05~\mathrm{ns/day}$ for PPPM, corresponding to an $\approx 2.9\times$ speedup. Thus, ESP samples the same NPT ensemble at substantially lower computational cost.

\subsection{Long-time simulations of a transmembrane bovine \texorpdfstring{\boldmath$bc_1$}{bc1} complex}
\begin{figure}[!htbp]
\centering
\includegraphics[width=0.92\linewidth]{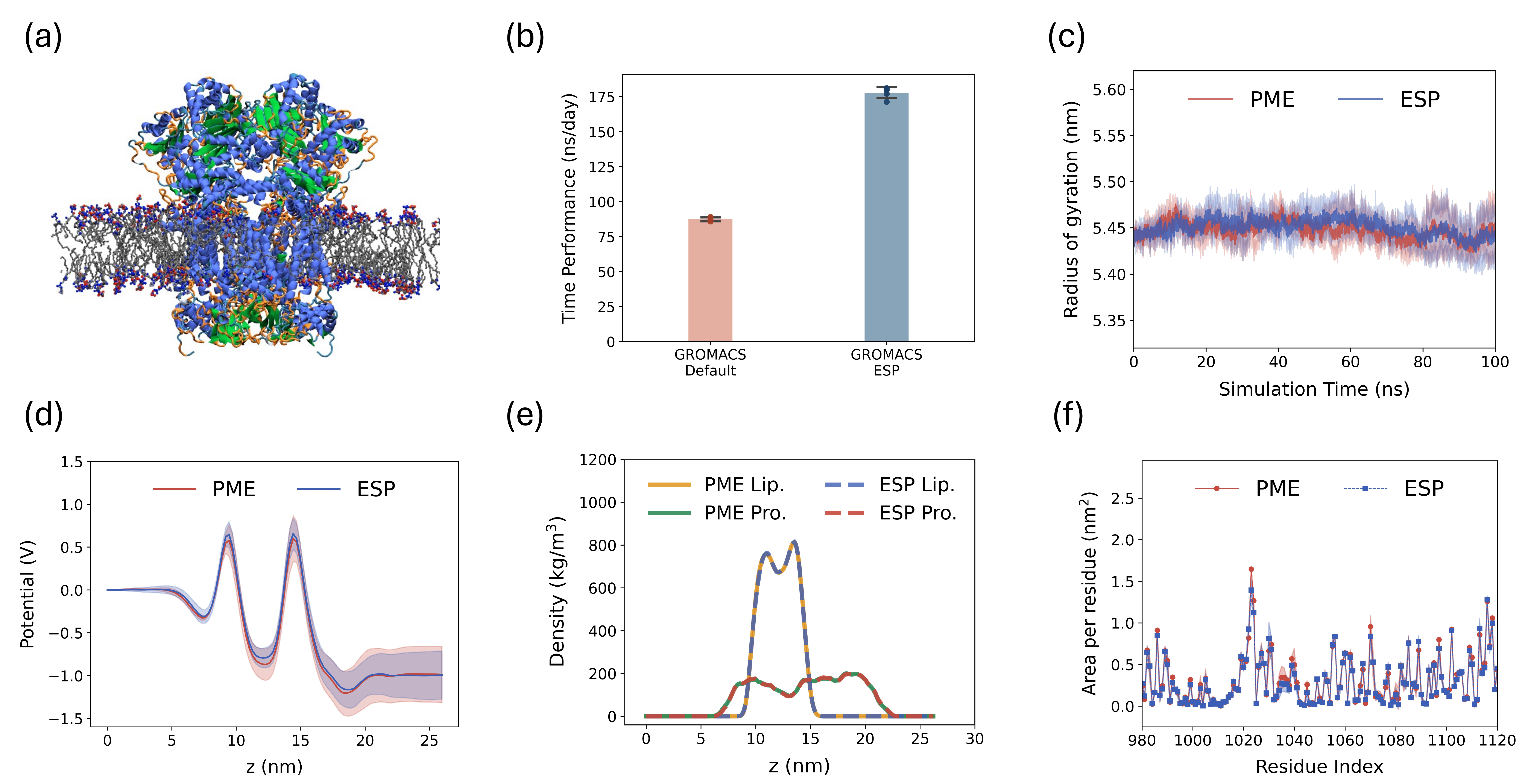}
\caption{\sf Comparison of NPT simulations for a transmembrane bovine cytochrome $bc_1$ complex ($809{,}997$ atoms). Results are from $100~\mathrm{ns}$ production runs (after equilibration) in GROMACS using PME and ESP. {\bf (a)} Snapshot of the system. Atoms on the lower and upper lipid surfaces are shown as deep blue and red spheres, respectively, based on the sign of their valence. Hydrophobic lipid tails are shown in gray (partially displayed for clarity). The protein is colored by secondary structure: $\alpha$-helices (blue), $\beta$-sheets (green), turns (yellow), and coils (cyan). {\bf (b)} Performance comparison across methods, including default and tuned GROMACS settings and ESP. {\bf (c)} Protein radius of gyration versus time; shaded bands indicate confidence intervals estimated from five independent runs. {\bf (d)} Electrostatic potential profile along the $z$-axis. {\bf (e)} Density profiles of lipids and protein along the $z$-axis. {\bf (f)} Solvent-accessible surface area (SASA) per residue for residues 980--1120 (transmembrane region).}
\label{fig:transmembrane}
\end{figure}

To test ESP in a complex biophysical setting, we simulate the bovine heart mitochondrial cytochrome $bc_1$ complex, a large transmembrane protein system (Fig.~\ref{fig:transmembrane}a). The cytochrome $bc_1$ complex is a key component of the electron transport chain, essential for ATP synthesis, and a drug target for antibiotics~\cite{xia1997crystal}. Each monomer contains 11 subunits, and the complex forms a homodimer in bovine mitochondria. We use the GROMOS 53A6 force field for proteins and lipids~\cite{oostenbrink2004biomolecular} and the SPC water model~\cite{berendsen1981interaction}. The initial configuration is taken from the MemProtMD database~\cite{newport2019memprotmd} using PDB entry 1sqq. The system contains $809{,}997$ atoms, including the transmembrane $bc_1$ protein, 1{,}050 DPPC lipids, 237{,}771 water molecules, 833 Na$^{+}$, and 832 Cl$^{-}$. After energy minimization, we equilibrate the system in the NPT ensemble at $323~\mathrm{K}$ and $1~\mathrm{bar}$ for $200~\mathrm{ns}$ using a Nos\'e--Hoover chain thermostat~\cite{martyna1992nose} and a stochastic cell-rescaling barostat~\cite{bernetti2020pressure} with semi-isotropic pressure coupling. We then perform a $100~\mathrm{ns}$ production run with time step $\Delta t=2~\mathrm{fs}$ on 960 CPU cores. The rectangular periodic box has dimensions $18.19\times 19.59\times 26.22~\mathrm{nm}$ in the $x$, $y$, and $z$ directions.

We compare ESP with the native GROMACS PME method~\cite{Darden1993JCP}. For both solvers, we truncate the short-range Coulomb and Lennard--Jones interactions at $r_c=0.9~\mathrm{nm}$, following the recommended force-field settings. For PME, we use the default parameters recommended in the GROMACS manual: spreading order $P=5$ and Fourier spacing $0.12~\mathrm{nm}$, which corresponds to a practical relative force error of approximately $4\times10^{-4}$. For ESP, we set $\Delta=4\times10^{-4}$, $P=5$, and Fourier spacing $0.26~\mathrm{nm}$ to match the same error level. Under these settings, ESP achieves approximately $180~\mathrm{ns/day}$, corresponding to an $\sim 2\times$ speedup over the default PME configuration (Fig.~\ref{fig:transmembrane}b).

For transmembrane proteins, accurate pressure control is important for maintaining stable secondary structure and lipid bilayers. To assess accuracy, we compute the protein radius of gyration, electrostatic potential profiles of non-solvent components along the $z$-axis, density profiles of lipids and protein along the $z$-axis, and solvent-accessible surface area (SASA) for residues 980--1120 (Figs.~\ref{fig:transmembrane}c--f). Across these metrics, ESP matches the native GROMACS results while using fewer than one-tenth of the Fourier-space grid points required by the default PME setup.

\section{Concluding remarks}\label{sec::ConcludingRemarks}
In this work, we develop an NPT-compatible formulation of Ewald summation with prolates (ESP) for charged systems under periodic boundary conditions. Building on our earlier ESP framework for spectrally accurate energy and force evaluation, the present work derives consistent pressure-tensor formulas for isotropic, semi-isotropic, anisotropic, and fully flexible cells, with explicit treatment of self-energy and the zero-frequency contribution. A key computational outcome is that the long-range pressure can be evaluated using a single forward FFT followed by diagonal scaling, whereas force evaluation requires both forward and inverse transforms.

Compared with current mesh-Ewald methods based on Gaussian splitting and B-spline spreading, ESP employs prolate spheroidal wave functions (PSWFs) for both splitting and spreading. This choice yields a more compact reciprocal-space representation at a prescribed accuracy, reducing both arithmetic work and all-to-all communication. Across bulk-water benchmarks, ESP exhibits spectral convergence and strong scaling with near-linear growth in problem size. Long-timescale simulations of a concentrated LiTFSI electrolyte and a transmembrane $bc_1$ complex further demonstrate that ESP reproduces structural, dynamical, and thermodynamic observables while using fewer than one-tenth of the Fourier grid points required by PME/PPPM in mainstream software, resulting in an overall $2$--$3\times$ speedup for NPT simulations. These results indicate that ESP is a practical and efficient option for large-scale NPT simulations on modern high-performance computing systems.

\section*{Acknowledgments}
The work of J. L. was partially supported by the National Natural Science Foundation of China (grant No. 12401570) and the China Postdoctoral Science Foundation (grant No. 2024M751948). The authors would like to thank the Scientific Computing Core at the Flatiron Institute for their support and for providing essential computational resources.
The Flatiron Institute is a division of the Simons Foundation.

\bibliographystyle{abbrv}
\bibliography{esp_npt}

\appendix

\section{Mathematical preliminaries}
\subsection{Basic properties of the PSWF function}\label{app::pswf}
Let $c>0$ be a real parameter. The prolate spheroidal wave function (PSWF) is an eigenfunction of the compact integral operator $\mathscr{F}_c:L^{2}[-1,1]\rightarrow L^2[-1,1]$ defined by
\begin{equation}\label{eq::formula}
\mathscr{F}_c[\varphi](x)=\int_{-1}^{1}\varphi(t)\,e^{icxt}\,dt.
\end{equation}
We denote the eigenvalues by $\lambda_0,\lambda_1,\ldots$ and order them so that $|\lambda_n|\ge |\lambda_{n+1}|$ for all $n\ge 0$. Let $\psi_n^{c}$ be an eigenfunction associated with $\lambda_n$, i.e.,
\begin{equation}\label{eq::psidefinit}
\lambda_n\psi_{n}^{c}(x)=\int_{-1}^{1}\psi_{n}^{c}(t)\,e^{icxt}\,dt,
\qquad x\in[-1,1],\; n\ge 0.
\end{equation}
We normalize $\psi_{n}^{c}$ so that $\|\psi_{n}^{c}\|_{L^2[-1,1]}=1$. It is known~\cite{osipov2013prolate} that $\{\psi_n^{c}\}_{n\ge 0}$ are real-valued, orthonormal, and complete in $L^2[-1,1]$.

A key feature of PSWFs is their joint time--frequency concentration. In particular, among all $L^2$ functions supported on $[-1,1]$ with unit $L^2$ norm, the \emph{order-zero} PSWF $\psi_0^{c}$ uniquely maximizes the fraction of Fourier energy contained in the band $[-c,c]$ (equivalently, it minimizes the $L^2$ energy outside $[-c,c]$); higher-order PSWFs provide subsequent maximizers subject to orthogonality constraints~\cite{osipov2013prolate,slepian1983sirev}. This optimal concentration property motivates using $\psi_0^{c}$ as a near-optimal compactly supported window.

Finally, PSWFs admit a simple Fourier relation on the band $[-c,c]$. If $\psi_n^{c}$ is extended by zero outside $[-1,1]$ and we consider its Fourier transform restricted to $|k|\le c$, then
\begin{equation}\label{eq::pswf_fourier}
\widehat{\psi}_{n}^{c}(k)=\lambda_n\,\psi_{n}^{c}(k/c),
\qquad |k|\le c.
\end{equation}
In particular, for $n=0$ this shows that the band-limited Fourier transform of the compactly supported window $\psi_0^{c}$ reproduces the same function (up to the scalar factor $\lambda_0$) under the scaling $k\mapsto k/c$, analogous to the self-reproducing property of Gaussians under the Fourier transform.

\subsection{Convolution theorem}\label{app:convolution}
Let $f(\bm{r})$ and $g(\bm{r})$ be periodic functions on $\Omega$ with Fourier coefficients $\widehat{f}(\bm{k})$ and $\widehat{g}(\bm{k})$, respectively, under the Fourier convention used in Section~\ref{subsec::pswfdecomp}. Define the (periodic) convolution
\begin{equation}
u(\bm{r}):=(f*g)(\bm{r})=\int_{\Omega} f(\bm{r}-\bm{r}')\,g(\bm{r}')\,d\bm{r}'.
\end{equation}
Then the Fourier coefficients satisfy
\begin{equation}
\widehat{u}(\bm{k})=\widehat{f}(\bm{k})\,\widehat{g}(\bm{k}).
\end{equation}

\subsection{Plancherel's theorem}\label{app:PlancherelTheorem}
Let $f(\bm{r})$ and $g(\bm{r})$ be periodic functions on $\Omega$ with Fourier coefficients $\widehat{f}(\bm{k})$ and $\widehat{g}(\bm{k})$, respectively. Under the same Fourier convention, Plancherel's identity reads
\begin{equation}
\int_{\Omega} f(\bm{r})\,\overline{g(\bm{r})}\,d\bm{r}
=\frac{1}{V}\sum_{\bm{k}} \widehat{f}(\bm{k})\,\overline{\widehat{g}(\bm{k})},
\end{equation}
where $V=\det(\bm{h})$ is the cell volume and the overline denotes complex conjugation. The special case $f\equiv g$ is the Parseval--Plancherel identity.

\subsection{Fourier transform of radially symmetric functions}\label{App:Fourier}
Assume $f(\bm{x})$ is integrable on $\mathbb{R}^3$ so that its Fourier transform exists. If $f$ is radially symmetric, $f(\bm{x})=f(|\bm{x}|)$, then $\widehat{f}$ is also radially symmetric and can be written as
\begin{equation}
\widehat{f}(k)=4\pi\int_{0}^{\infty} f(r)\,\frac{\sin(kr)}{kr}\,r^{2}\,dr,
\qquad k=|\bm{k}|.
\end{equation}

\section{Distribution functions in the NPT ensemble}\label{app::DistributionFunction}
We briefly review distribution functions for the NPT ensemble under several common parameterizations of the simulation cell. Detailed derivations can be found in~\cite{frenkel2001understanding,bond2007molecular}.

For a system coupled to a heat bath and a pressure reservoir, the NPT weight of a microstate with energy $E$ and volume $V$ is proportional to $\exp[-\beta(E+PV)]$. For fixed cell shape, the isothermal--isobaric partition function is~\cite{frenkel2001understanding}
\begin{equation}
\Delta(N,P,T)
=\int_{0}^{\infty} e^{-\beta PV}\,Q(N,V,T)\,dV
=\int e^{-\beta(E+PV)}\,d\bm{r}_{\mathrm{tot}}\,d\bm{p}_{\mathrm{tot}}\,dV,
\end{equation}
where $Q(N,V,T)$ is the canonical (NVT) partition function at volume $V$.

When the cell shape is allowed to fluctuate, one must specify an integration measure for the shape variables. We write the cell tensor as $\bm{h}=V^{1/3}\bm{h}_0$ with $\det(\bm{h}_0)=1$, so that $\bm{h}_0$ parameterizes shape at fixed volume. Typical choices are as follows.

\begin{enumerate}
\item \textbf{Semi-isotropic coupling.}
Let $\Omega$ be a rectangular box with square base area $A$ (coupled isotropically in $x$ and $y$) and height $L$ (along $z$), so that $V=AL$. This setup is widely used for membrane simulations. Introducing reduced shape variables
$A_0:=V^{-2/3}A$ and $L_0:=V^{-1/3}L$ with the constraint $A_0L_0=1$, a convenient shape measure is
$dA_0\,dL_0\,\delta(A_0L_0-1)$, where $\delta(\cdot)$ is the one-dimensional Dirac delta. The partition function can be written as
\begin{equation}
  \Delta(N,P,T)=\int_{0}^{\infty} e^{-\beta PV}\,Q_{\mathrm{semi}}(N,V,T)\,dV,
\end{equation}
where
\begin{equation}
Q_{\mathrm{semi}}(N,V,T)=\int_{0}^{\infty}\!\!\int_{0}^{\infty} Q(N,V,T)\,\delta(A_0L_0-1)\,dA_0\,dL_0.
\end{equation}
Equivalently, in variables $(A,L)$ one may write
\begin{equation}
\Delta(N,P,T)=\int e^{-\beta(E+PAL)}\,d\bm{r}_{\mathrm{tot}}\,d\bm{p}_{\mathrm{tot}}\,dA\,dL,
\end{equation}
where $E$ may depend on the cell shape through $(A,L)$. If a nonzero surface tension $\gamma_0$ is imposed (the $NP\gamma_0T$ ensemble~\cite{zhang1995computer}), the weight becomes $\exp[-\beta(E+PAL-\gamma_0A)]$.

\item \textbf{Anisotropic coupling in a rectangular cell.}
For a rectangular cell with side lengths $(L_x,L_y,L_z)$, the partition function in variables $\{\bm{r}_{\mathrm{tot}},\bm{p}_{\mathrm{tot}},L_x,L_y,L_z\}$ can be written as
\begin{equation}
\Delta(N,P,T)=\int e^{-\beta\!\left(E+P L_xL_yL_z\right)}
\,d\bm{r}_{\mathrm{tot}}\,d\bm{p}_{\mathrm{tot}}\,dL_x\,dL_y\,dL_z.
\end{equation}

\item \textbf{Fully flexible cell.}
For a general cell tensor $\bm{h}$ with $\det(\bm{h})>0$, a commonly used form is
\begin{equation}
\Delta(N,P,T)=\int_{\det(\bm{h})>0} \det(\bm{h})^{-2}\,e^{-\beta P \det(\bm{h})}\,Q(N,V,T)\,d\bm{h}.
\end{equation}
The factor $\det(\bm{h})^{-2}$ leads to an additional term in the extended Hamiltonian used to sample the correct ensemble; it arises from the choice of measure in $\bm{h}$-space.
\end{enumerate}

To sample the NPT ensemble, the thermostat and barostat must generate the corresponding target distribution. In our GROMACS implementation, we combine ESP with the stochastic cell-rescaling barostat~\cite{bernetti2020pressure}. For example, in the isotropic case the barostat evolves the log-volume variable $\varepsilon$ according to
\begin{equation}
d\varepsilon
=-\frac{\beta_T}{\tau_P}\left(P_0-P_{\mathrm{ins}}\right)\,dt
+\sqrt{\frac{2k_{\mathrm{B}}T\,\beta_T}{V\,\tau_P}}\,dW,
\end{equation}
where $\varepsilon=\log(V/V_0)$ with reference volume $V_0$, $P_0$ is the target external pressure, $\tau_P$ is the barostat time constant, $\beta_T$ is the isothermal compressibility, and $W$ is a standard Wiener process. In our LAMMPS implementation, we combine ESP with the MTK barostat~\cite{martyna1994constant}. More generally, ESP can be used with other standard thermostat and barostat schemes.

\section{Unified treatment of force and pressure calculations}\label{app::force}
In NPT simulations, time integration requires both the forces and the instantaneous pressure tensor. Although PSWF-based force evaluation was introduced in our earlier ESP work~\cite{liang2025acceleratingfastewaldsummation}, here we emphasize that force and pressure can be computed within a single particle--mesh workflow, reusing intermediate quantities and FFTs.

From Eqs.~\eqref{eq::splitting}--\eqref{eq::selfen}, the force on particle $i$ is the negative gradient of the total electrostatic energy,
\begin{equation}\label{eq::Fc}
\begin{split}
\bm{F}(\bm{r}_i)
&=\sum_{\bm{n}}{}^{\prime}\sum_{j=1}^{N}q_iq_j\,
\mathcal{F}_{\mathcal{N}}\!\left(|\bm{r}_{ij}+\bm{h}\bm{n}|\right)
\frac{\bm{r}_{ij}+\bm{h}\bm{n}}{|\bm{r}_{ij}+\bm{h}\bm{n}|^{3}}
-\sum_{\bm{k}\neq\bm{0}}\frac{q_i}{V}\,\widehat{\mathcal{F}}^{c}(\bm{k})\,
\bm{k}\,\Im\!\left[e^{-i\bm{k}\cdot\bm{r}_i}\rho(\bm{k})\right] \\
&=: \bm{F}_{\mathcal{N}}(\bm{r}_i)+\bm{F}_{\mathcal{F}}(\bm{r}_i),
\end{split}
\end{equation}
where $\mathcal{F}_{\mathcal{N}}$ is defined in Eq.~\eqref{eq::mathFN}, $V=\det(\bm{h})$, and $\Im(\cdot)$ denotes the imaginary part. Here, $\bm{F}_{\mathcal{N}}$ and $\bm{F}_{\mathcal{F}}$ denote the short-range and long-range components evaluated in real space and Fourier space, respectively.

\noindent \textbf{Short-range force and pressure.}
Combining Eqs.~\eqref{eq::Fc} and~\eqref{eq::PinsN} shows that the short-range force and short-range pressure share the same pairwise interaction kernel. Let $\mathcal{N}_i$ be the neighbor list of particle $i$ within cutoff radius $r_c$. Define the pairwise short-range force
\begin{equation}
\bm{F}_{\mathcal{N},ij}
:=q_iq_j\,\mathcal{F}_{\mathcal{N}}(|\bm{r}_{ij}|)\,\frac{\bm{r}_{ij}}{|\bm{r}_{ij}|^{3}}.
\end{equation}
Then
\begin{equation}
\bm{F}_{\mathcal{N}}(\bm{r}_i)=\sum_{j\in\mathcal{N}_i}\bm{F}_{\mathcal{N},ij},
\qquad
\bm{P}_{\mathrm{ins},\mathcal{N}}
=-\frac{1}{2V}\sum_{\bm{n}}{}^{\prime}\sum_{i,j=1}^{N}
\bm{F}_{\mathcal{N},ij}^{(\bm{n})}\otimes(\bm{r}_{ij}+\bm{h}\bm{n}),
\end{equation}
where $\bm{F}_{\mathcal{N},ij}^{(\bm{n})}$ denotes the corresponding interaction evaluated with the periodic shift $\bm{h}\bm{n}$. In practice, the same neighbor-list loop can accumulate both $\bm{F}_{\mathcal{N}}$ and $\bm{P}_{\mathrm{ins},\mathcal{N}}$ at essentially no extra cost.

\noindent \textbf{Long-range force: two lattice differentiation options.}
As in standard mesh-Ewald methods~\cite{Darden1993JCP,essmann1995smooth}, there are two spectrally accurate approaches for differentiating the long-range contribution.

\begin{enumerate}
\item \textbf{$i\bm{k}$-differentiation.}
Differentiating in Fourier space amounts to multiplying each Fourier coefficient by $-i\bm{k}$. In this case, the long-range force can be expressed as
\begin{equation}
\bm{F}_{\mathcal{F}}(\bm{r}_i)
=-\frac{q_i}{V}\sum_{\bm{k}\neq\bm{0}}\widehat{W}(\bm{k})\,e^{i\bm{k}\cdot\bm{r}_i}\,
\widehat{\bm{\rho}}_{\mathrm{diag}}^{\,\mathrm{force}}(\bm{k}),
\end{equation}
where
\begin{equation}\label{eq::rhodiagforce}
\widehat{\bm{\rho}}_{\mathrm{diag}}^{\,\mathrm{force}}(\bm{k})
=i\bm{k}\,\widehat{\mathcal{F}}^{c}(\bm{k})\,|\widehat{W}(\bm{k})|^{-2}\,\widehat{\rho}_{\mathrm{grid}}(\bm{k}).
\end{equation}
Starting from $\widehat{\rho}_{\mathrm{grid}}(\bm{k})$ (obtained after spreading and the forward FFT), Eq.~\eqref{eq::rhodiagforce} is applied mode-by-mode (diagonal scaling), and three inverse FFTs are then used to recover the three components of the vector field in real space. The force on particles is obtained by a standard gathering step.

\item \textbf{Analytical differentiation.}
Alternatively, one can differentiate the real-space representation obtained after applying the inverse FFT. Using Eq.~\eqref{eq::UF}, the long-range energy can be written as
\begin{equation}\label{eq::UFAnal}
\begin{split}
U_{\mathcal{F}}
&=\frac{1}{2V}\sum_{i=1}^{N}q_i\sum_{\bm{k}\neq\bm{0}}
\widehat{W}(\bm{k})\,e^{i\bm{k}\cdot\bm{r}_i}\,
\widehat{\rho}_{\mathrm{diag}}^{\,\mathrm{energy}}(\bm{k})
=\frac{1}{2}\sum_{i=1}^{N}q_i\int_{\Omega}\rho_{\mathrm{diag}}^{\,\mathrm{energy}}(\bm{r})\,W_{*}(\bm{r}_i-\bm{r})\,d\bm{r},
\end{split}
\end{equation}
where
\begin{equation}
\widehat{\rho}_{\mathrm{diag}}^{\,\mathrm{energy}}(\bm{k})
=\widehat{\mathcal{F}}^{c}(\bm{k})\,|\widehat{W}(\bm{k})|^{-2}\,\widehat{\rho}_{\mathrm{grid}}(\bm{k}).
\end{equation}
Differentiating the real-space form in Eq.~\eqref{eq::UFAnal} yields
\begin{equation}\label{eq::AD}
\bm{F}_{\mathcal{F}}(\bm{r}_i)
=-\frac{q_i}{2}\int_{\Omega}\rho_{\mathrm{diag}}^{\,\mathrm{energy}}(\bm{r})\,\nabla_{\bm{r}_i} W_{*}(\bm{r}_i-\bm{r})\,d\bm{r}.
\end{equation}
In practice, $\rho_{\mathrm{diag}}^{\,\mathrm{energy}}$ is obtained using one inverse FFT, and the gathering step evaluates $\nabla_{\bm{r}_i} W_{*}(\bm{r}_i-\bm{r})$ at particle positions in the three coordinate directions.
\end{enumerate}

LAMMPS~\cite{plimpton1995fast} typically uses $i\bm{k}$-differentiation, whereas GROMACS~\cite{berendsen1995gromacs} typically uses analytical differentiation. The two approaches have different conservation properties: $i\bm{k}$-differentiation conserves momentum but not energy, whereas analytical differentiation conserves energy (for sufficiently small time steps) but not momentum~\cite{Hockney1988Computer}. In analytical differentiation, conservation of center-of-mass momentum can be enforced by an additional correction; this may introduce a small energy drift, which can be mitigated by subtracting the self-force~\cite{ballenegger2012convert}. In the NPT setting considered here, these effects are typically negligible because thermostats and barostats already regulate energy and volume fluctuations. For efficiency and consistency with our ESP workflow, we use analytical differentiation in this work.

\end{document}